\documentclass[DIV=12]{scrartcl}
\usepackage{xcolor}

\newcommand{\blue}[1]{{\color{blue!80!black} \noindent #1}}

\RequirePackage[algo2e,algoruled,norelsize, linesnumbered, lined,commentsnumbered]{algorithm2e}

\usepackage{amsmath,amssymb,amsthm}
\usepackage{booktabs}
\usepackage{caption}
\usepackage{subcaption}
\usepackage{lipsum}
\usepackage{amsfonts}
\usepackage{graphicx}
\usepackage{epstopdf}
\usepackage{amsopn,numprint}
\usepackage{placeins}
\usepackage[hidelinks]{hyperref}

\newcommand{\TR}[1]{{#1}}

\newcommand{\algow}{0.2mm}
\hypersetup{
    colorlinks,
    linkcolor={red!50!black},
    citecolor={blue!50!black},
    urlcolor={blue!80!black}
}

\ifpdf
  \DeclareGraphicsExtensions{.eps,.pdf,.png,.jpg}
\else
  \DeclareGraphicsExtensions{.eps}
\fi
\graphicspath{{.}{figures/}}

\usepackage{float}
\usepackage[within=section]{newfloat}
\DeclareFloatingEnvironment{algorithm}

\usepackage{listings}
\usepackage{units}

\def\Re{\text{Re}}
\def\bu{\boldsymbol u}
\def\bv{\boldsymbol v}
\def\bb{\boldsymbol b}

\def\bv{\boldsymbol v}
\def\bx{\boldsymbol x}
\def\bV{\boldsymbol V}
\def\bf{\boldsymbol f}
\def\bchi{\boldsymbol\chi}
\def\bsigma{\boldsymbol\sigma}
\def\beps{\boldsymbol\epsilon}

\def\div{\mbox{div}\,}

\usepackage[twoside]{fancyhdr}

\title{An adaptive finite element multigrid solver using GPU acceleration}
 
\author{M. Liebchen\thanks{Institute of Analysis and Numerics, OVGU Magdeburg, Germany (\url{manuel.liebchen@ovgu.de}).}
\and R. Jendersie\thanks{Institute of Analysis and Numerics, OVGU Magdeburg, Germany
    (\url{robert.jendersie@ovgu.de}).}
  \and U. Kaya\thanks{Institute of Analysis and Numerics, OVGU Magdeburg, Germany
    (\url{utku.kaya@ovgu.de}).}
  \and C. Lessig\thanks{European Centre for Medium-Range Weather Forecasts
    (\url{christian.lessig@ecmwf.int}).}
  \and T. Richter\thanks{Institute of Analysis and Numerics, OVGU Magdeburg, Germany
    (\url{thomas.richter@ovgu.de}).}
  }

\begin{document}

\maketitle

\pagestyle{fancy}

\begin{abstract}
Adaptive finite elements combined with geometric multigrid solvers are one of the most efficient numerical methods for problems such as the instationary Navier-Stokes equations.
Yet despite their efficiency, computations remain expensive and the simulation of, for example, complex flow problems can take many hours or days. 
GPUs provide an interesting avenue to speed up the calculations due to their very large theoretical peak performance. 
However, the large degree of parallelism and non-standard API make the use of GPUs in scientific computing challenging.
In this work, we develop a GPU acceleration for the adaptive finite element library Gascoigne and study its effectiveness for different systems of partial differential equations.
Through the systematic formulation of all computations as linear algebra operations, we can employ GPU-accelerated linear algebra libraries, which simplifies the implementation and ensures the maintainability of the code while achieving very efficient GPU utilizations.
Our results for a transport-diffusion equation, linear elasticity, and the instationary Navier-Stokes equations show substantial speedups of up to 20X compared to multi-core CPU implementations. 
% with, for example, a speedup of more than 35X for the Navier-Stokes equations. 
\end{abstract}

%adaptive finite elements, GPU parallelization, Navier-Stokes equations

%%%
\section{Introduction}
\label{sec:introduction} 

The combination of adaptive finite elements and geometric multigrid solvers is one of the most efficient approaches for the numerical approximation of partial differential equations. Adaptive mesh refinement schemes based on numerical a posteriori error estimates allow for optimal refinement~\cite{Stevenson2005} while geometric multigrid solvers provide the solution to the arising linear systems of equations in linear complexity~\cite{Hackbusch1985}. 
For many partial differential equations, the combination yields an overall approach with optimal computational complexity for a prescribed accuracy.

The mathematical and computational optimality come at the price of complex algorithms and a lack of regular structure in the computations. 
This is due to the adaptivity of the discretization, resulting in an irregular mesh and stencils, and the multi-scale nature of the solvers. 
Furthermore, memory is often accessed in an unstructured way, hampering the parallelization of the computations.

These challenges are one of the central reasons that parallel accelerator co-processors, such as GPUs and TPUs, are still only rarely used in numerical simulations, including adaptive multigrid methods. 
Accelerators are highly attractive since the computational power of CPUs has stagnated in the last decade while those accelerators has increased dramatically. 
% Indeed, today a large part of the compute power is typically found in the accelerator units. 
% These are still hardly used in numerical simulation, including for adaptive multi-grid methods, due to the difficulties in exploiting the highly parallel accelerators efficiently, especially in complex computations.
% where CPUs and accelerators have to work together to attain optimal performance.
% Typical benchmark problems reach only small fractions of the theoretical peak performance. 

Accelerators use task- and data-parallelism as well as specialization to achieve a very large peak performance. 
While a high degree of parallelism is inherent in some problems, e.g. in computer graphics and deep learning, it poses challenges for the implementation of adaptive algorithms in numerical linear algebra and scientific computing.
One direction to address this gap is by changing to algorithms with highly regular computations, such as lattice Boltzmann methods in fluid mechanics~\cite{OBRECHT2013252,Ribbrock2010}.

\TR{The performance characteristics of modern accelerator hardware, in particular GPUs, fundamentally differ from those of traditional CPUs. While GPUs offer massive parallelism and high arithmetic throughput, the performance in computations is in practice is often limited by bandwidth and latency between GPU-RAM and computational units, as well as by the cost of data transfers between host and device. These constraints pose significant challenges for classical finite element methods, which typically rely on sparse matrix storage and irregular memory access patterns to be efficient. As a consequence, a substantial body of recent work has focused on algorithmic reformulations of finite element methods that reduce memory traffic and increase arithmetic intensity to make more efficient use of GPUs.
A key development in this direction is the use of matrix-free methods, where discrete operators are applied without explicitly assembling and storing sparse matrices. Instead, operator entries are recomputed on the fly, typically using tensor-product structure and sum-factorization techniques for higher-order elements. An  overview of such approaches is given in~\cite{10.5555/3108096.3108097}, where the potential of tensor-product evaluations and sum factorization for high-performance finite element computations is discussed, albeit not yet in the context of multigrid solvers. These ideas have since become central to GPU-oriented finite element design due to their favorable compute-to-memory ratio~\cite{Munch2023}.

On GPUs, particularly high performance for multigrid methods has been demonstrated for fixed-stencil discretizations, which can be interpreted as a limit case of matrix-free computation. Thereby the application of an operator reduces to a small number of regular, local stencil operations with predictable memory access patterns. Recent studies such as~\cite{Antepara2024,Richardson2025} report excellent GPU efficiency for such methods, including for example  $p$-multigrid techniques.  While these results represent an important performance benchmark, their applicability to general, adaptive finite element methods is limited by the strong reliance on regular grid structures and fixed discretizations. 

Also general finite element frameworks have made significant progress in leveraging GPUs. The MFEM library~\cite{Andrej2024} represents a modern, GPU-centric approach that supports matrix-free higher-order finite elements and provides $p$-multigrid preconditioning for Krylov subspace solvers. A distinguishing feature of MFEM is its partial realization of an end-to-end GPU workflow, in which most stages of the solver—including operator application and multigrid components—can be executed on the device, thereby minimizing data transfers. A similar direction is pursued in the deal.II library~\cite{Munch2023}. There GPU support has been steadily expanded within a matrix-free high-order finite element framework. Recent developments emphasize performance portability and the integration of accelerator backends into existing multigrid and matrix-free infrastructures. While these efforts demonstrate the feasibility of combining geometric multigrid with GPUs in general-purpose finite element libraries, they also highlight the complexity of designing algorithms that remain efficient across different hardware architectures and use cases.

A central difficulty in matrix-free multigrid methods on GPUs lies in the design of effective smoothers. Due to the lack of an explicitly assembled matrix, smoothers are often restricted to point-wise operations such as Jacobi or Chebyshev iterations. This may lead to reduced multigrid efficiency, especially for high-order discretizations~\cite{Munch2023}. Notable progress is in recent contributions on multigrid methods for high-order discontinuous Galerkin discretizations~\cite{Cui2025,Cui2025a}, where efficient GPU-based vertex-patch smoothers are developed. These approaches demonstrate that carefully designed local solvers exploiting tensor-product structure can overcome some of the traditional limitations of matrix-free smoothing. The same concepts have also been successfully applied to the Stokes equations, using multigrid directly in the coupled velocity–pressure~\cite{Cui2025b}. In~\cite{Heuveline2012}, the authors investigate the GPU parallelization of several smoothers within a matrix-based geometric multigrid method. The implementation allows for locally refined unstructured meshes, but, the obtained speedup is limited and does not substantially exceed that of multi-core CPUs

In contrast to matrix-free methods, matrix-based approaches on GPUs are less common, primarily due to the high memory footprint of sparse matrices and the irregular memory access patterns of sparse matrix–vector products. Nevertheless, these approaches remain attractive because they can often be integrated more easily into existing simulation frameworks. In~\cite{Xu2025}, a two-level multigrid scheme with ILU-based smoothing on GPUs is investigated, demonstrating that meaningful acceleration can still be achieved for certain problem classes. A broader perspective on GPU usage in large-scale solvers is provided in the context of the PETSc framework~\cite{Mills2025}. The authors discuss solver design choices and scaling behavior toward exascale systems.

The primary advantage of matrix-free methods lies in their ability to drastically reduce memory usage by recomputing operator entries when needed. In situations where matrix-free solvers are not readily available—such as in complex multiphysics applications or legacy codes—mixed-precision techniques offer an alternative strategy. A general discussion of mixed-precision algorithms on GPUs, CPUs, and hybrid systems is given in~\cite{Abdelfattah2021}. In~\cite{RudaTurek} mixed-precision is made possible by using a hierarchical approach to modify the linear systems, denoted as pre-handling. In the context of multigrid solvers, mixed precision has been explored particularly for algebraic multigrid methods in~\cite{Tsai2023}. The potential of these approaches stems not only from reduced memory bandwidth requirements, but also from the significantly higher memory and compute throughput of low-precision arithmetic on modern GPUs.

For algebraic multigrid methods, GPU accelerated libraries exist, e.g. AmgX~\cite{Naumov2015} and Ginkgo~\cite{ginkgo-toms-2022}. NVIDIA~\cite{NVIDIAMG} has also presented a highly efficient geometric multigrid solver for 3d linear elliptic problems. A highly optimized pressure Poisson solver in GPUs is discussed in~\cite{GMRESAMG2024}.

An alternative memory-reduction strategy orthogonal to matrix-free methods is proposed in~\cite{Boukhris2023}. There sparse matrices are represented in a block-wise manner using embedded stencil-like structures. This approach aims to combine some of the flexibility of matrix-based methods with the reduced memory footprint and regular access patterns typically associated with stencil computations.

Finally, a fundamentally different direction is the use of machine learning techniques, in particular deep neural networks, as surrogate models or solver components in scientific simulations. Such models can be executed very efficiently on accelerator hardware and can typically exploit from reduced numerical precision, leading to speedups of several orders of magnitude in specific applications, as demonstrated, for example, in~\cite{Bi:2023aa,Margenberg2024}. While these approaches are conceptually distinct from classical multigrid methods, they are increasingly explored as complementary tools in large-scale simulation workflows.

A minimally invasive use of GPUs in classical finite element simulation tools is also possible but often yields little to no acceleration~\cite{5191718,Goddeke}. The obtained speedup was moderate but the work describes how a transparent implementation into an existing software can be achieved. The consequences of mixed precision arithmetic were also discussed in this work. 
In~\cite{GEVELER2013327}, the same authors presented a geometric multigrid method that is based entirely on sparse matrix-vector multiplications and can thus be easily and efficiently implemented on different hardware using suitable libraries. Applied to stationary linear differential equations, the solver exhibited up to a fivefold speedup relative to CPU-based systems.
}

In this article, we show that highly efficient adaptive finite elements with multigrid solvers can be implemented efficiently on GPU accelerators without limiting the flexibility of the method and without the large implementation effort that GPU-implementations can involve.
For this, we use a matrix-based formulation that allows us to use the cuBLAS and cuSPARSE libraries~\cite{cuSPARSE} that are highly optimized while the code remains close to those with BLAS calls on the CPU.
Hence, the CPU and GPU versions can retain the same structure and differ only in the linear algebra function calls and additional CPU-to-GPU memory transfer in the GPU version.
It, nonetheless, allows for the use of custom CUDA kernels when necessary and we demonstrate that the use of simple native CUDA code can provide significant speedups for operations that do not naturally map to (sparse) linear algebra.
Due to these optimizations, the largest performance bottleneck comes from the memory transfer between CPU and GPU, in particular, if system matrices are assembled on the CPU and then transferred to the GPU. In Section~\ref{sec:results} we study numerical examples, where the repeated assembly can be avoided such that matrices are copied only once for all time steps. This reduces the overhead and the remaining time for memory transfer is negligible.

In our final code, large parts of the computations have been transferred to the GPU and thus the slowdown through the transfer is also small.

We demonstrate the generality and flexibility of our approach by applying it to two linear elliptic problems, namely the transport-diffusion equation in 2d and a 3d linear elasticity problem. The discretization of these two equations results directly in a linear system of equations that can be approximated with the GPU-accelerated geometric multigrid solver. 
As a third example, we consider the instationary Navier-Stokes equations. For these, we first derive an explicit pressure-correction method, which can be represented by matrix-vector multiplications and a pressure-Poisson problem, the latter one being solved with the multigrid solver. Special attention is given to the nonlinearity of the Navier-Stokes equations. Through a reformulation, also this will be approximated by a matrix-vector multiplication with a pre-computed sparse matrix on the GPU.
For all three examples, we obtain significant speed-ups between 5X and 20X for our final GPU parallelizations.
At the same time, the systematic use of (sparse) linear algebra leads to easily understandable and maintainable code for the GPU computations.

\paragraph{Outline}
The remainder of the article is structured as follows. 
In Sec.~\ref{sec:adaptive_fem} we introduce the mathematical notation and briefly present the finite element discretization. There, we also describe the geometric multigrid method. 
In Section~\ref{sec:accel}, we will discuss accelerators such as GPUs and describe the specifics of the hardware and how they need to be reflected in the algorithms. We also introduce software libraries that facilitate the use of accelerator hardware. 
Section~\ref{sec:cuda_multigrid} describes our implementation of the multigrid process based on the cuSPARSE library~\cite{cuSPARSE}. Finally, in Section~\ref{sec:results}, we present numerical test problems and discuss the results.
We conclude in Section~\ref{sec:conclusions}.

\paragraph{Main contributions}
We describe the integration of GPU acceleration into the Gascoigne 3d general-purpose finite element library. As the software did not consider accelerator hardware in its design, the focus is on making the most efficient use of the hardware while minimising interventions into the existing implementation. For this reason, we base our GPU acceleation on cuSPARSE and require custom kernels only for a small set of operations that cannot be expressed efficiently using linear algebra. 
Additionally, we present interpolation-based finite element techniques that enable us to formulate certain nonlinear terms with matrix-vector products based on sparse matrices with small stencils. This avoids numerical quadrature and extends the applicability of the cuSPARSE-based implementation to a larger class of problems. Our insight for the GPU acceleration is not specific to Gascoigne 3d and will hence likely generalize to other finite element libraries. 

\section{Adaptive Finite Element Discretization and Multigrid Solver}
\label{sec:adaptive_fem}

\subsection{Finite element discretization}

We denote a domain by $\Omega\subset\mathbb{R}^d$, where $d=2$ or $d=3$ is the dimension. By $\Omega_h$ we denote a finite element mesh consisting of $N_h$ quadrilaterals or generalized (allowing curved faces) hexahedras. The elements $T\in\Omega_h$ all arise from a reference element $\hat T$
\[
\hat T_T:\hat T\mapsto T
\]
that is the unit square in 2d and the unit cube in 3d. We will consider isoparametric finite elements, where the mapping itself comes from the finite element space. Let
\[
Q^r:=\{ x_1^{\alpha_1}\cdots x_d^{\alpha_d},\; \alpha_i\in\mathbb{N},\; 0\le \alpha_i\le r\}
\]
be the space of polynomials of maximal degree $r$ in each coordinate. Then, $\hat T_T\in [Q^r]^d$ and we define the finite element spaces of degree $r$ as
\[
V_h^{(r)} := \{ \phi\in C(\bar\Omega)\,:\,
\phi\circ T_T^{-1} \in Q^r\quad\forall T\in \Omega\}.
\]
We assume shape regularity in the sense that $\|\nabla T_T\|\cdot \|\nabla T_T^{-1}\|\le c$ uniform in $h>0$, see~\cite[Sec. 4.2]{Richter2017}. Structural regularity is relaxed by allowing local mesh refinement with at most 1 hanging node per face, again, see~\cite[Sec 4.2]{Richter2017} for the specific realization in the finite element library Gascoigne 3d~\cite{Gascoigne}.

Systems of partial differential equations such as the Navier-Stokes equations or elasticity will be discretized with equal-order finite elements, i.e. $\mathbf{V}_h = [V_h^{(r)}]^c$, where $c\in\mathbb{N}$ is the number of components. Gascoigne 3d combines these components locally. Taking the 3d Navier-Stokes equations as example this means that the vector is represented as a matrix with entries $x_{i,c}$ where $i$ refers to the grid node and $c$ to the component ($c=1$ pressure, $c=2,3,4$ velocities). The system matrix also has double indexing. Each matrix entry $A_{ij}$ is itself a (dense) matrix $A_{ij}\in\mathbb{R}^{n_c\times n_c}$ with $n_c=4$ in the case of 3d Navier-Stokes. Likewise, each vector entry $v_i$ is itself a vector $v_i\in\mathbb{R}^{n_c}$.
Keeping the solution components together is beneficial in terms of cache efficiency for problems where one obtains large system matrices like Navier-Stokes ($n_c=4$) or 3d elasticity, where $n_c=6$ when the velocity and deformation fields are combined~\cite{BraackRichter2006d,P1}. 
The memory layout of this approach is not easily transferred to GPUs as standard libraries like cuSPARSE do not support such blocked matrices. Section~\ref{sec:gpu:mem} will give details. 

\subsection{Adaptivity and geometric multigrid hierarchy}\label{ada}

\begin{algorithm}[t]
  \small
  \caption{Geometric multigrid solver}
  \label{algo:gmg}
\rule{\textwidth}{\algow}
Given a hierarchy of meshes $\Omega_l$ and corresponding finite element spaces $\mathbf{V}_l$ for $l=0,\dots,L$ and an initial value $x_L^0\in \mathbf{V}_L$. For $n=1,\dots,N_{max}$ iterate
  \[
  x_L^{(n)} = GMG(L,x_L^{(n-1)},b_L)
  \]
  where $GMG(l,x_l,b_l)$ is recursively defined as:
\[
\begin{aligned}
&\text{\textbf{Step 0:} Coarse mesh problem}
&&&\text{if }l=0&\text{ return }A_0^{-1}b_0\\
&\text{\textbf{Step 1:} Pre-smooth}
&&&x_l' &= {\cal S}_l(x_l,b_l)\\
&\text{\textbf{Step 2:} Restrict residual}
&&&d_{l-1} &= {\cal R}_{l-1}(b_l-A_l x_l')\\
&\text{\textbf{Step 3:} Recursive coarse mesh correction}
&&&y_{l-1} &= GMG(l-1,0,d_{l-1})\\
&\text{\textbf{Step 4:} Prolongate update}
&&&x_l'' &= x_l' + {\cal P}_{l-1}y_{l-1}\\
&\text{\textbf{Step 5:} Post-smooth}
&&& &\text{return }{\cal S}_l(x_l'',b_l)
\end{aligned}
\]
%\rule{\textwidth}{\algow}
\end{algorithm}

Adaptivity of the computations is realized by means of hierarchical local mesh refinement. 
If an element $T\in\Omega_h$ is chosen for refinement, it is replaced by 4 elements in 2d and 8 elements in 3d. We allow a level jump of 1 for neighboring elements. The actual refinement can therefore extend further into the domain. Faces and edges on elements with a level jump therefore have inner nodes that are unknowns only on one side. These nodes are called hanging nodes and are replaced by interpolations of their direct neighbours in the finite element approach. In the original Gascoigne, hanging node data used for this interpolation is stored in a sparse structure only covering exactly the hanging nodes. As preparation for GPU parallelization, the data is now represented in a sparse matrix. This matrix $\mathbf{H}_h$ is the diagonal identity in all rows belonging to standard degrees of freedom, whereas in rows belonging to hanging nodes, the matrix contains the interpolation weights such that the matrix vector product $\mathbf{y} = \mathbf{H}_h \mathbf{x}$ performs all interpolations at once. The computational effort is $O(N)$ ($N$ being the number of unknowns), regardless if hanging nodes are present. However, the constant is small and the overhead in relation to a treatment of only the hanging nodes is negligible. The transpose $\mathbf{H}_h^T$ plays an important role in the assembly of the matrix and variational residuals. If, for example, an integral $f_i=(f,\phi_i)$ is calculated for $i=1,\dots,N$, all test functions $\phi_i$, including those in hanging nodes, are processed first. The result vector $\mathbf{f}_h=(f_i)_i$ is then multiplied by $\mathbf{H}_h^T$ so that the test functions belonging to hanging nodes are correctly taken into account.  

Mesh refinement can be based either on a priori knowledge or on a posteriori error estimates, the latter usually in the context of the dual weighted residual method~\cite{BeckerRannacher2001,BraackRichter2006d}.

\begin{figure}[t]
\begin{center}
    \includegraphics[width=\textwidth]{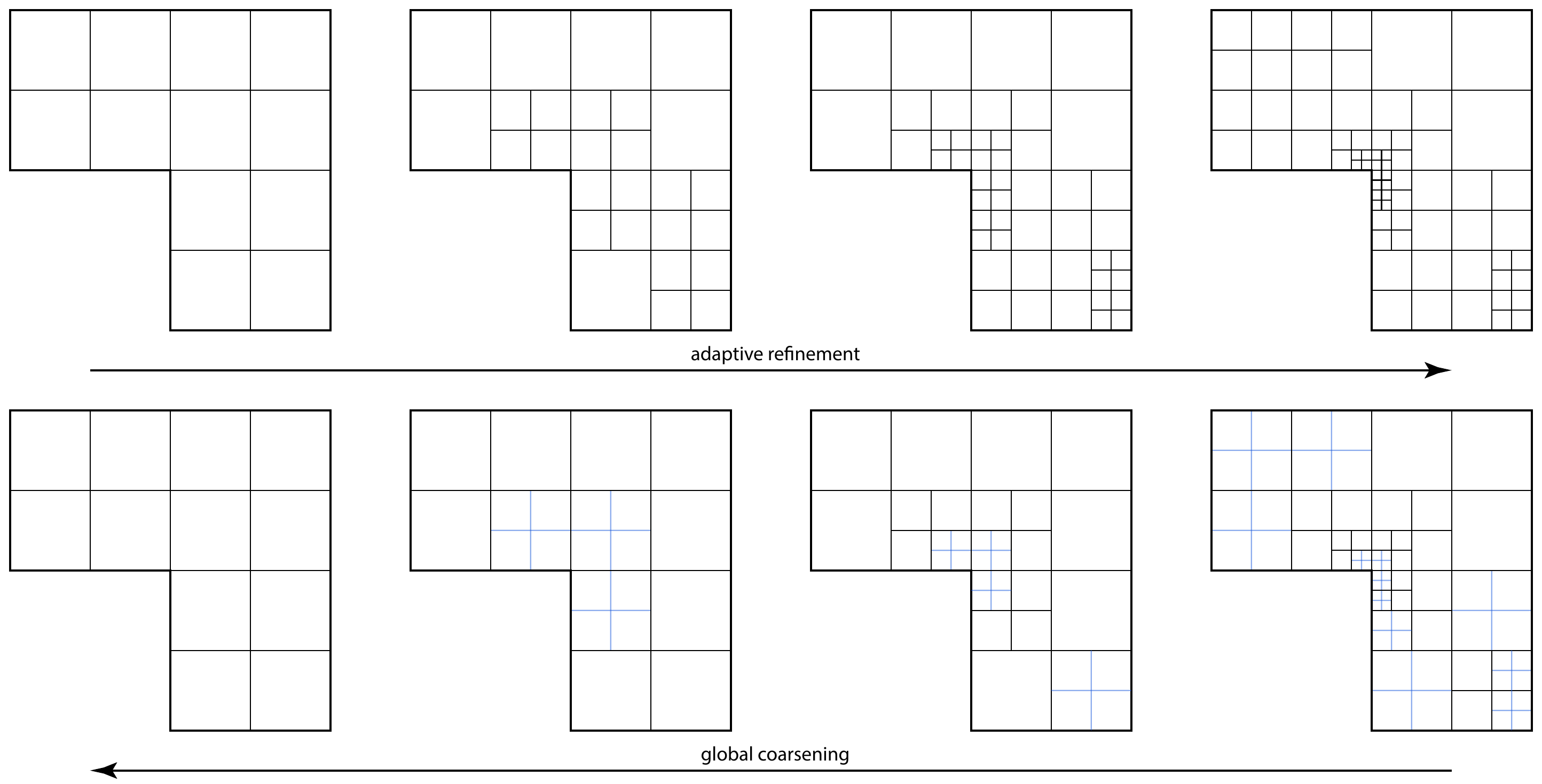}
\end{center}
\caption{Top row: hierarchy of adaptively refined 2:1 meshes with hanging nodes. Bottom row: global coarsening resulting in the multigrid hierarchy. In each step as many refinements are taken back as possible (shown in blue) for a rapid decrease of the mesh complexity.}
\label{fig:coarsening}
\end{figure}

Successive mesh refinement results in a hierarchy of meshes. Instead of using this hierarchy for the setup of the geometric multigrid solver, we start with the finest mesh $\Omega_h=:\Omega_L$ and recursively coarsen it until we reach a coarse mesh $\Omega_H=:\Omega_0$ with small complexity $N_H$. The mesh $\Omega_H$ is either the coarse starting mesh or any finer mesh chosen to give an optimal balance between recursive multigrid performance and fast direct solution of a small problem. By $L$ we denote the resulting number of multigrid layers. In each coarsening step, as many refinements as possible are taken back, yielding a hierarchy $\Omega_h=\Omega_L \succeq \Omega_{L-1}\succeq\dots \Omega_0=\Omega_H$. 
To be precise, an element is coarsened, if it belongs to a group of four (in 3d eight) elements on the same mesh level that all arise from splitting the same common father element. Details are given in~\cite{BeckerBraack2000a}.
We follow a global multigrid approach where each intermediate level $\Omega_l$ spans the complete domain, see~\cite{BeckerBraack2000a}. The main advantage of the global multigrid is its simplicity. An approximation of the entire problem can be created on each grid level. This allows, for example, the introduction of global constraints, e.g. a normalisation of the pressure, $\int_\Omega p\,\text{d}x=0$, in the context of incompressible flows or compliance with material balances $\sum_i c_i=1$ for chemical reactions, see~\cite{P2}. The disadvantage of the global approach is the greater effort required, as the multigrid smoother covers the entire grid on each level.  However, if one analyses the effort in terms of the accuracy that can be achieved with adaptive grids and not in terms of the number of degrees of freedom, it becomes clear that the global approach has no significant disadvantage in practical applications~\cite{BeckerBraack2000a,J1,BraackRichter2006d}.

Fig.~\ref{fig:coarsening} provides a sketch of the adaptive mesh refinement and the global coarsening procedure that result in two different hierarchies of meshes.

The multigrid solver is the standard $V$-cycle. Coarse mesh problems are either solved exactly or approximated using a couple of smoothing steps. 
The geometric multigrid iteration is either used as linear solver or as pre-conditioner in a GMRES solver~\cite{Saad1996}. The GMRES solver is more robust and usually required for problems that are transport dominated or when meshes with large element aspect ratios are used.

%%%%%%%%%%%%%%%%%%%%%%%%%%%%%%%%%%%%%%%%%%%%%%%%%%%%%%%%%%%%%%%%%%%%%%%%%%%%%%%%%%%%%%
\section{Data-parallel accelerators}\label{sec:accel}

Accelerator co-processors, such as GPUs from NVIDIA and AMD and TPUs from Google, are used in a wide range of applications, such as computer gaming and neural networks training. Due to their very large compute power, they also play an increasingly important role in scientific computing. 
The large compute power is achieved through a very high degree of parallelization as well as specialization. 
For example, NVIDIA's GPUs combine MIMD parallelism with up to 132 streaming multi-processors with data-parallelism with a logical, pipelined width of $1024$ and a hardware width of $32$ (and $8$ on the latest hardware).
Additionally, significant compute power is provided by dedicated matrix-matrix multiplications engines known as tensor cores and introduced for neural network training.
Accelerator co-processors also differ from conventional processors through native hardware support of lower precision data types such as half precision and even 8-bit computations on the latest generation of chips.
Accelerator co-processors as used in this work come with their own, dedicated memory hierarchy starting from a RAM, and with typically two layers of caches.
To perform computations on the accelerator, data hence needs to first be transferred to the accelerator RAM through a comparatively slow memory interface. 
% The transfer of input data to the accelerator and of results back to the CPU RAM can hence easily be a bottleneck in computations. 

The programming of accelerator hardware is typically challenging due to the high degree of parallelism and narrow ``fast paths'' on the hardware (e.g. because of smaller caches compared to CPUs) that are difficult to optimize for and can change from hardware generation to hardware generation. This is compounded by the separate memory hierarchy on the accelerator and the slow interface to CPU RAM, that can easily become a bottleneck.
Furthermore, native software libraries, such as CUDA and ROCm, are usually vendor specific and not standardized, with often also an incomplete documentation. 
This makes their use by non-experts challenging and leads to a high maintenance effort when always the latest features should be used.

Different libraries have been proposed to aid with the development of GPU-accelerated software.
OpenACC and OpenMP provide pragma-based access to the compute power of GPUs that is very simple to use but also limited in the potential for optimizations. 
SYCL~\cite{syclmain} and Kokkos~\cite{Trott2022} provide a middle-ground with more complex usage but also more flexibility.
Machine learning libraries such as PyTorch~\cite{paszke2017automatic} and jax~\cite{jax2018github} follow a new paradigm for GPU programming where the application code is specified in a high-level language and a backend, such as XLA~\cite{Sabne2020_XLA} or ATen~\cite{paszke2017automatic}, generates highly efficient, hardware-specific code.
This concept is sometimes also referred to as domain-specific language, e.g.~\cite{Lawrence2018}.
Recently, Triton~\cite{Tillet2019} was proposed as a complement to existing machine learning libraries.
It still allows for highly simplified programming of accelerators compared to CUDA and ROCm but provides more control than PyTorch and jax, and hence typically yields more efficient code.

\LinesNumbered%
\SetFuncSty{textbf} \SetCommentSty{textsf} \SetKwInOut{Input}{Input} %
\SetKwProg{Function}{function}{}{end} %
\SetKwFunction{MG}{Multigrid}%
\SetKwFunction{rhs}{Rhs}%

For specific applications, such as linear algebra or discrete Fourier transform, also highly optimized libraries exist. These provide very high performance while, through the targeted use case, also allow for much simpler usage than, e.g., CUDA.
In our work, we will build on a linear algebra formulation of our adaptive PDE solvers.
This allows us to use the cuBLAS~\cite{cuBLAS} and cuSPARSE libraries~\cite{cuSPARSE} that provide high performance while being simple to use with a software interface that is roughly comparable to those of the conventional BLAS library.
Furthermore, cuBLAS and cuSPARSE can still be combined with native CUDA code when necessary, i.e. when an operation cannot be expressed efficiently in linear algebra.

\section{CUDA multigrid}
\label{sec:cuda_multigrid}

In this section we discuss our extension of the adaptive finite element library Gascoigne 3d to use GPU accelerators for a wide range of computations. 
The extensions have been designed with two use cases in mind: first, non-stationary but linear systems of partial differential equations, where a single system matrix is reused in every time step or at least in many time steps, and, second,  nonlinear problems that can be formulated in terms of fixed matrices that do not need re-assembly.
The second case includes a pressure-projection based solvers for the nonstationary nonlinear Navier-Stokes equations but can be extended to more complex flow problems such as the Boussinesq approximation~\cite{PAMMBoussinesq}.
Both problem classes have in common that matrices are assembled only once and the complete workflow can be formulated in terms of matrix vector products and the solution of linear systems of equations. 
Algorithms~\ref{algo:linear} and~\ref{algo:navierstokes} show the typical workflow for the problem types. Computations highlighted in blue are completely performed on the GPU, orange marks transfer of data between GPU and CPU.

\begin{algorithm}[t]
    \caption{Nonstationary Linear Problem   \label{algo:linear}}
\rule{\textwidth}{\algow}
    Basic Initialization\tcp*{\texttt{init}}
    Assemble Matrices\tcp*{\texttt{init}}
    {\color{orange!60!black}CopyToGPU(Matrix)}\tcp*{\color{black}\texttt{copy}}
    \For{time-iter}{
    {\color{blue}Assemble Right Hand Side}\tcp*{\texttt{rhs}}
%    {\color{orange!60!black}CopyToGPU(RighHandSide)}\tcp*{\color{black}\texttt{copy}}
    {\color{blue}SolveLinearProblem}\tcp*{\texttt{solve}}}
    {\color{orange!60!black}CopyFromGPU(Solution)}\tcp*{\color{black}\texttt{copy}}
% \rule{\textwidth}{\algow}
\end{algorithm}

It is possible to also port the computations of general nonlinear problems to a GPU. However, keeping the flexibility and support of different discretizations and triangulations is not as easily handled for the ill-structured setup of typical finite element assemblies of system matrices and residuals in the general nonlinear case. 
We hence defer this to future work.

\begin{algorithm}[t]
    \label{algo:driven_cavity}
    \caption{Explicit Navier-Stokes Pressure Correction\label{algo:navierstokes}}
\rule{\textwidth}{\algow}
    Basic Initialization\\
    Assemble Matrices\\
    {\color{orange!60!black}CopyToGPU(Matrix)}\tcp*{\color{black}\texttt{copy}}
    {\color{blue}\For{time-iter}{
        Assemble Momentum right Hand Side{\color{black}\tcp*{\texttt{mom-rhs}}}
        Explicit Solve Momentum Problem{\color{black}\tcp*{\texttt{mom-solve}}}
        Assemble Pressure Right Hand Side{\color{black}\tcp*{\texttt{pres-rhs}}}
        \For{GMRES iteration}{
        Multigrid Pressure 
        {\color{black}\tcp*{\texttt{pres-solve}}}}
        Update Pressure{\color{black}\tcp*{\texttt{pres-up}}}
    Correct Velocity
    }}
    {\color{orange!60!black}CopyFromGPU(Solution)}\tcp*{\color{black}\texttt{copy}}
% \rule{\textwidth}{\algow}
\end{algorithm}

\subsection{General architecture}

Our principal approach for the CUDA parallelization of Gascoigne 3d is to formulate all operations in terms of (potentially sparse) matrix-vector and matrix-matrix products and to realize these with the highly optimized cuBLAS~\cite{cuBLAS} and cuSPARSE~\cite{cuSPARSE} libraries.
This approach differs from the existing implementation in Gascoigne 3d where many operations such as the handling of hanging nodes or the multigrid mesh transfer are based on the local connectivity of the degrees of freedom as function of mesh and discretization. However, this functionality can easily be realized as matrix-vector products, including for the handling of hanging nodes and all mesh transfer operations.
This is, e.g., the approach taken from the outset in the deal.II finite element library~\cite{dealII}. 
% Hence, some computations in Gascoigne 3d needed to be formulated in matrix-vector form first. 

All high level operations in Gascoigne 3d are performed on abstract interface classes to vectors and matrices to allow for MPI parallelization~\cite{P1,P2}.
This has been retained in the CUDA-version so that the high-level control flow and complex algorithms (such as Newton, multigrid, GMRES) are unchanged since they only operate on the interfaces. 
The approach minimizes the changes that were overall required and allows to flexibly retain both the CPU and the GPU backend.

% . This allows to retain most of Gascoigne 3d's software framework. Only those classes must be complemented with CUDA versions where direct memory access is required and where data either resides on the GPU or on the CPU. 
% Other code implementing high-level algorithms is unchanged.  

%TODO - Wir machen ja die hängenden Knoten und den Mehrgittertransfer per Matrix-Vektor Produkt. Auch hier kannst Du vielleicht schon einmal ein paar Stichpunkte einfügen.

% Our CUDA extension 

% A multigrid solver has the advantage that for some calculations not the full resolution is needed computations on coarser levels are computationally more efficient. 
% The mesh hierarchy is thereby generated by subdividing a coarsest mesh. We only use topogically rectangular 2D structures or cube like 3D structures where a regular subdivision is trivial. 
% For the transformation of the {\color{red} vectors} between grid sizes, non-symmetric BSR matrices {\color{red}[cite]} are used that prolongate or project the values according to the subdivision structure.

% {\color{red} GPU?}

% All operations written as matrix-vector for CUDA, in contrast to Gascoigne.
% Simplify implementation and exploit optimizations. 

%%%%%%%%%%%%%%%%%%%%%%
\subsubsection{Storage}\label{sec:gpu:mem}

\begin{figure}[t]
  \begin{center}
    \includegraphics[width=\textwidth]{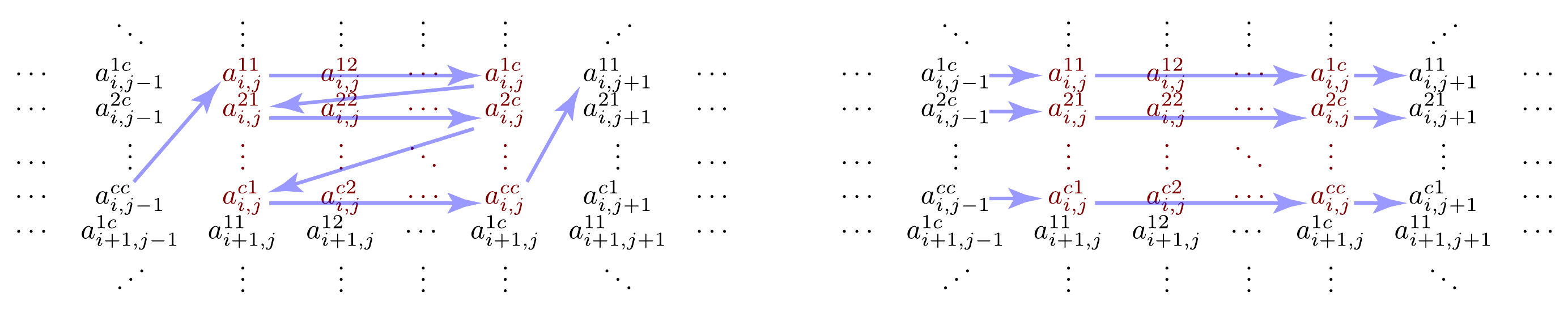}
  \end{center}
  \caption{Storage of matrix on CPU (left) and GPU (right). On the CPU we employ a block-wise ordering of the matrix entries clustering the $n_c\times n_c$ components of a system of PDEs. On the GPU standard CSR format is used. }
  \label{fig:matrixmem}
\end{figure}

Gascoigne 3d uses a block-wise memory concept in which the various components of a partial differential equation (e.g. the 3 deformation unknowns in 3d solid mechanics or the pressure and the three velocity components in fluid mechanics) are stored next to each other. Matrices are double-indexed and the outer index $(i,j)$ refers to the mesh node yielding the entry $a_{i,j}\in\mathbb{R}^{n_c\times n_c}$ that stores (aligned in memory) the local couplings between the solution components as a matrix. 
In general, this setup helps to efficiently use caching~\cite{BraackRichter2006d,P1}. Using cuSPARSE this block-matrix format is not supported. While vectors are directly transferred to the GPU using the same memory layout, the inner matrix blocks must be resolved and matrices are stored in the usual CSR format. Fig.~\ref{fig:matrixmem} shows a visualization of the matrix storage concept on the CPU (left) and the GPU (right).

The representation of a matrix in GPU memory is directly allocated when the matrix is created and assembled on the CPU. This helps to better exploit the asynchronous data transfer. Our implementation is directed at problems where the matrix stays fixed such that this transfer is required just once. For vectors, the data transfer takes place when a \texttt{cudaActivate()} function is called on the object owing the data, e.g. a class steering the multigrid solver. 
Data is transferred back to the CPU when the corresponding \texttt{cudaDeactivate()} function is called.

% The operations are thereby only performed when necessary, i.e. when the respective object not already resides on the target device for the operation.
% It also ensures that an appropriate error message is generated when not all arguments for an operation are on the same device.

%TODO - Das ganze Environment, also CPU / GPU / Speicher, etc. beschreiben. Das brauchen wir auf jeden Fall.

% Storage is an important topic in GPU-Acceleration. The time it takes to copy data from the CPU to the GPU or back is often large than those required for computation. As on the CPU, allocating memory also takes some time on the GPU. To keep the time of copying and reserving to a minimum,  the data being used is stored. In order to make it clear to the user of the library there is an explicit operation for switching to the GPU and coping particular data. This operation exists in both directions. 

%%%%%%%%%%
\subsubsection{Mapping (sparse) linear algebra to the GPU}

The matrix-matrix and matrix-vector multiplications that are the essential building blocks of CUDA-Gascoigne 3d are realized with cuBLAS and cuSPARSE, which provide an interface in close analogy to BLAS. 
For example, the function \texttt{cuSPARSESpMV} performs for a sparse version of BLAS's gemv, i.e. it computes $\alpha \, \mathrm{op}(A) x + \beta y$ where $\alpha$, $\beta$ are scalar, $x$, $y$ are vectors and $A$ is a matrix.
The approach reduces the implementation effort for our GPU-parallelization substantially and ensures it is maintainable, i.e. further developments of Gascoigne 3d can be implemented without substantially efforts also for the GPU version. 
% The matrices are thereby stored in 
% Also supported are block-sparse computations using the Blocked-Ellpack (ELL) format.

%%%%%%%%%%
\subsubsection{Custom CUDA kernels}

Almost all operations required for the adaptive finite element solvers of Gascoigne 3d can be expressed efficiently using linear algebra. 
However, we found some exceptions and performing these calculations on the CPU incurred a very high penalty due to the required data transfer between GPU and CPU.
We therefore implemented small parts of the computations in native CUDA, which avoided the extra memory transfers and hence the penalty.
This is facilitated by cuBLAS and cuSPARSE operating on raw CUDA pointers in device memory which can directly also be used in native CUDA.
Details of the implementation are given in Sec.~\ref{sec:ns}.  

The flexibility to combine cuBLAS, cuSPARSE and native CUDA is, in our opinion, an important feature to avoid CPU-GPU data transfer and to be able to flexibly implement algorithms with CUDA-Gascoigne-3d. The total number of required custom kernels is very small and serves as guideline for measuring the effort to port further applications to the GPU.

%%%%%%%%%%
\subsubsection{Geometric Multigrid CUDA}

As a concrete example, we describe in the following how geometric multigrid algorithm in Algo.~\ref{algo:gmg} is realized in CUDA-Gascoigne-3d. 
All computations are entirely performed on the GPU. Hence, data transfer overhead is incurred only at the beginning and end of the computations, see also Algorithms~\ref{algo:linear} and~\ref{algo:navierstokes}.  

All steps of the algorithm involve only elementary operations that can directly be formulated in cuSPARSE or cuBLAS. Mostly, matrix-vector products must be computed, e.g. for prolongation and restriction, but also for the smoother. To avoid memory transfer between CPU and GPU it is essential that the smoother can be performed completely on the GPU. At the moment we limit ourselves to very simple smoothers of Jacobi or block-Jacobi type that can be written as
\[
x^{(l+1)} = x^{(l)} + \omega S\big(b-Ax^{(l)}\big),
\]
where $\omega\in\mathbb{R}$ is a damping factor and $S\in\mathbb{R}^{n\times n}$ the smoothing operator, written as fixed matrix. More complex smoothers can, in principle, be implemented using native CUDA but we leave this to future work.

Mesh transfer operations are usually done locally. In the hierarchical setup of finite elements on quadrilateral or hexahedral meshes, prolongation of a solution to the next finer mesh is the usual embedding. To illustrate this, let $K$ be an element of the coarse mesh $\Omega_l$ and $K_1,\dots,K_p$ be the resulting fine elements on level $\Omega_{l+1}$, where $p=2^d$ with $d$ being the spatial dimension. Considering finite elements of degree $r$, $(r+1)^d$ unknowns are involved on mesh level $l$ and $(2r+1)^d$ on level $l+1$. Prolongation is then by means of 
\begin{equation}\label{prolongation}
u_i^{(l+1)} = \sum_{j=1}^{(r+1)^d} \chi_{ij}u_j^{(l)},\quad i=1,\dots,(2r+1)^d.
\end{equation}
The $(r+1)^d\times (2r+1)^d$ coefficients $\chi_{ij}$ are the same for each mesh element. \eqref{prolongation} can be written as one global matrix-vector product
\[
u_h^{(l+1)} = P_l u_h^{(l)}
\]
with $P_l\in\mathbb{R}^{N_{l+1}\times N_l}$ and its entries given by $\chi_{ij}$. 
With the number of distinct elements in $P_l$ being very small, the use of the matrix formulation would be sub-optimal on the CPU. On the GPU, however, it is performant and allows us to use cuSPARSE. 
The multigrid restriction $R_l:\Omega_{l+1}\to\Omega_l$ is the transpose of the prolongation, i.e. $R_l=P_l^T$ and hence also implemented using a sparse matrix in cuSPARSE.
Similarly to the mesh transfer, also for hanging nodes we formulate the operations for averaging and distributing the values of the solution vector as sparse matrix-vector products, see Section~\ref{ada}. 
% Likewise, these matrices are composed of repetitions some few different lines only, but this approach again reduced the need to introduce special code for the GPU version. 

% Here some context is missing

The coarse mesh problem in \textbf{Step 0} of Algorithm~\ref{algo:gmg} is not solved directly. Instead, we simply apply several steps of the smoothing iteration such that no additional infrastructure for GPU parallelization is required.

\subsubsection{GMRES as linear solver on the GPU}
\label{sec:gmres}

For reasons of numerical stability, e.g. on non-uniform meshes, the multigrid iteration can often not be used directly as linear solver. They are then employed as preconditioner in a GMRES iteration. Our GMRES solver follows the approach described in \cite[Section 6.5.3]{Saad1996}, using the modified Gram Schmidt algorithm for orthogonalization and Givens rotation for solving the resulting overdetermined linear system. Due to the multigrid preconditioning, we avoid stability issues in the Gram Schmidt orthogonalization and never have to perform more than 5-10 GMRES steps. The Gram Schmidt iteration is also completely run on the GPU using cuSPARSE and cuBLAS operations.
As we use a highly efficient geometric multigrid solver as preconditioner, a maximum of $n_{G}\le 10$ GMRES steps is usually required. Hence, the resulting overdetermined system for finding the GMRES solution is very small and its step-wise transformation to a diagonal matrix with Givens rotations is performed on the CPU. The cost for transferring these $n_G^2$ entries between GPU and CPU is small. Details on GMRES are given in~\cite{Saad1996}.

\section{Results}
\label{sec:results}

In this section, we present several applications where we employ the GPU-accelerated multigrid method presented in the last section. We start with two linear elliptic problems, namely the transport-diffusion equation in 2d and a 3d linear elasticity problem. The discretization of these two equations results directly in a linear system of equations that can be approximated with the geometric multigrid solver and falls into the class of systems handled by Algorithm~\ref{algo:linear}. 

As a third example, we consider the Navier-Stokes equations. For these, we first derive an explicit pressure-correction method, which can be represented by matrix-vector multiplications and a pressure-Poisson problem, the latter one being solved with the multigrid solver. Special attention is given to the nonlinearity of the Navier-Stokes equations. Through a reformulation, also this will be approximated by a product with a pre-computed sparse matrix and can easily be performed on the GPU.

\begin{figure}
    \centering
    \includegraphics[width=0.6\textwidth]{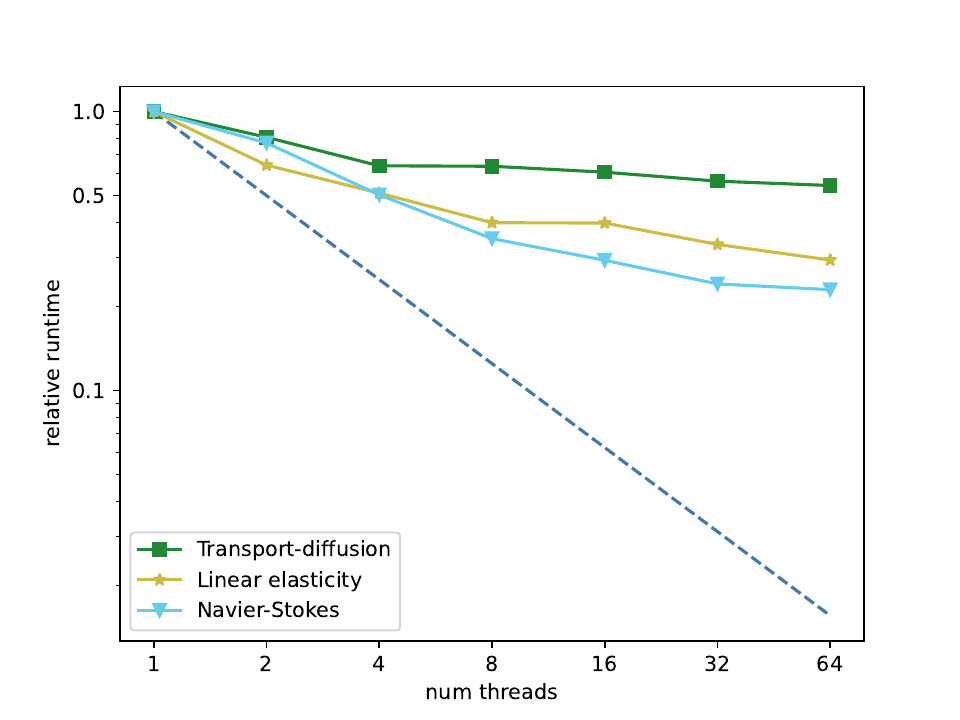}
    \caption{Strong scaling of the time-stepping parts of the largest problems described in Sections~\ref{sec:transport_diffusion}, \ref{sec:elasticitynum} and \ref{sec:ns}, with up to $2,000,000$ degrees of freedom on an AMD EPYC 7773X CPU.}
    \label{fig:cpu-scaling}
\end{figure}
% \begin{remark}[Measurement of computing time and hardware infrastructure]
%If not noted otherwise, all CPU computations are carried out on an Intel Xeon E5-2640 running at 2.40 Ghz using 8 threads. The GPU computations are performed on the same machine using an Nvidia V100 card with 16 GB of GPU RAM. When using the GPU to accelerate the multigrid solver, those parts of Gascoigne that reside on the CPU still run on 8 parallel threads.
All CPU computations are carried out on an AMD EPYC 7773X. For the relatively small problem sizes we consider, our scaling test in Fig.~\ref{fig:cpu-scaling} shows that there is only a minor benefit from using a large number of threads. The calculations are dominated by sparse matrix-vector multiplications and the efficiency is therefore limited by the memory bandwidth. This has already been shown in a detailed parallel analysis of fluid-structure interactions based on Gascoigne~\cite{Failer2020}. A transition to matrix-free finite elements—see, e.g., \cite{Kronbichler2019}—would thus offer the greatest potential for improving parallel efficiency. However, this would require a fundamental redesign of the current implementation.  We therefore limit the OpenMP based parallel code to 8 threads for the comparisons. The GPU computations are performed on the same machine using an NVIDIA H100 PCIe card with 80 GB of VRAM. When using the GPU to accelerate the multigrid solver, those parts of Gascoigne that reside on the CPU still run on 8 parallel threads.
We always indicate wall-clock times and separately specify the different contributions of Algorithms~\ref{algo:linear} and~\ref{algo:navierstokes}. 

\subsection{Transport-diffusion equation}
\label{sec:transport_diffusion}
\begin{figure}[t]
\centering
\includegraphics[width=.35\linewidth]{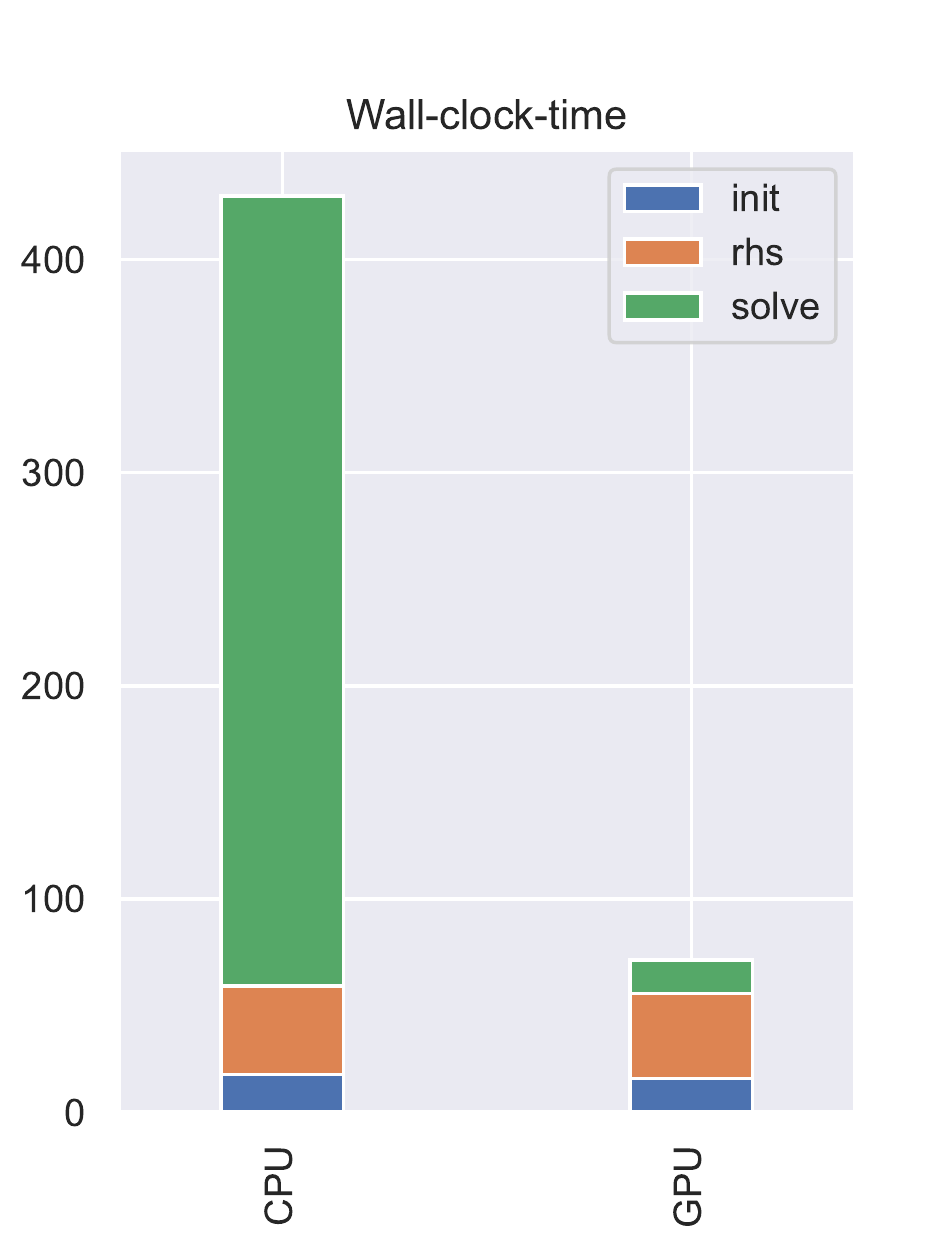}
  \includegraphics[width=.35\linewidth]{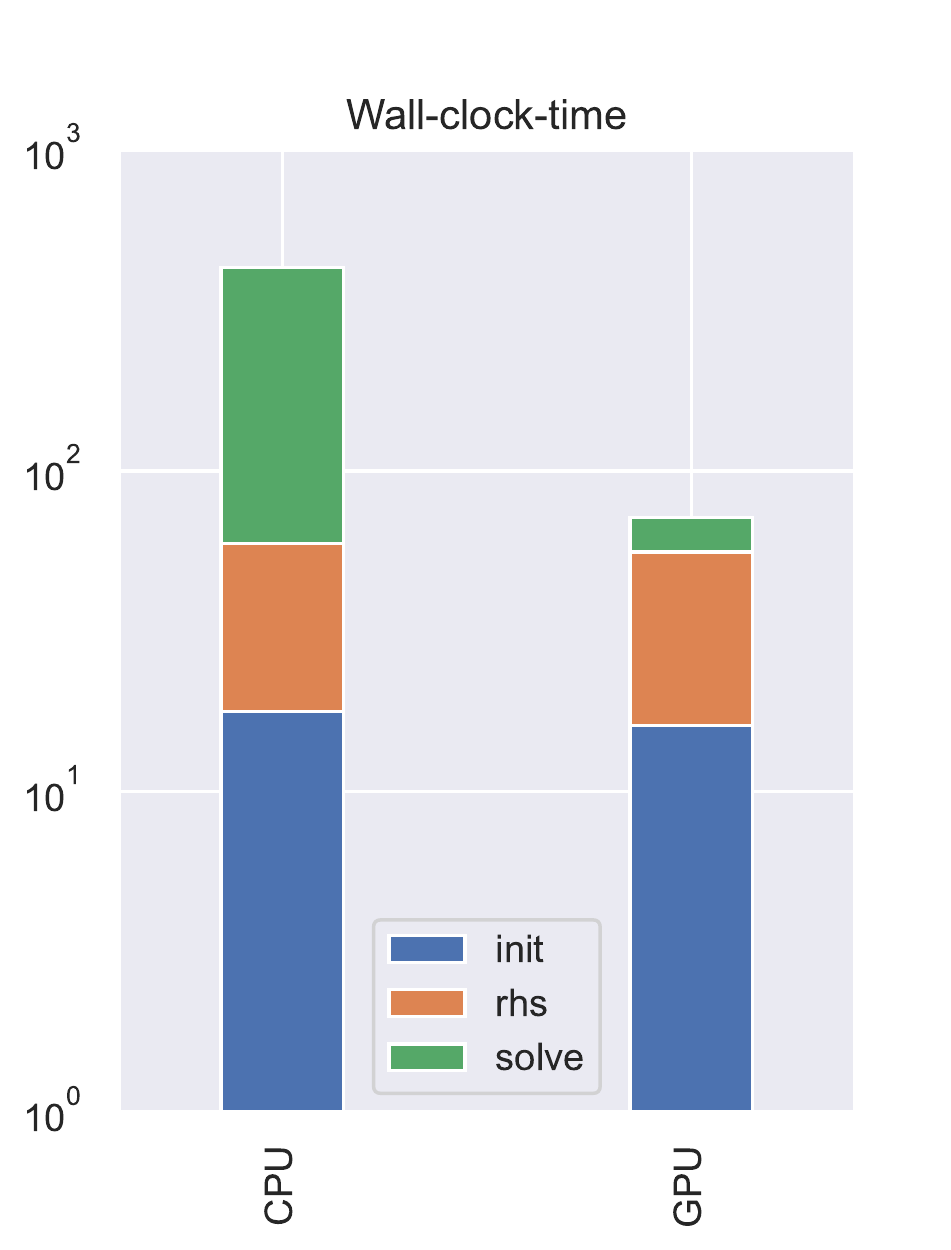}
  \caption{Transport-diffusion equation: Wall-clock times CPU vs. GPU on the finest mesh level with $1\,000\,000$ degrees of freedom (left: linear scale, right: log scale).}
  \label{fig:conv-dif-1}
\end{figure}

\begin{figure}[t]
 \includegraphics[width=0.48\linewidth]{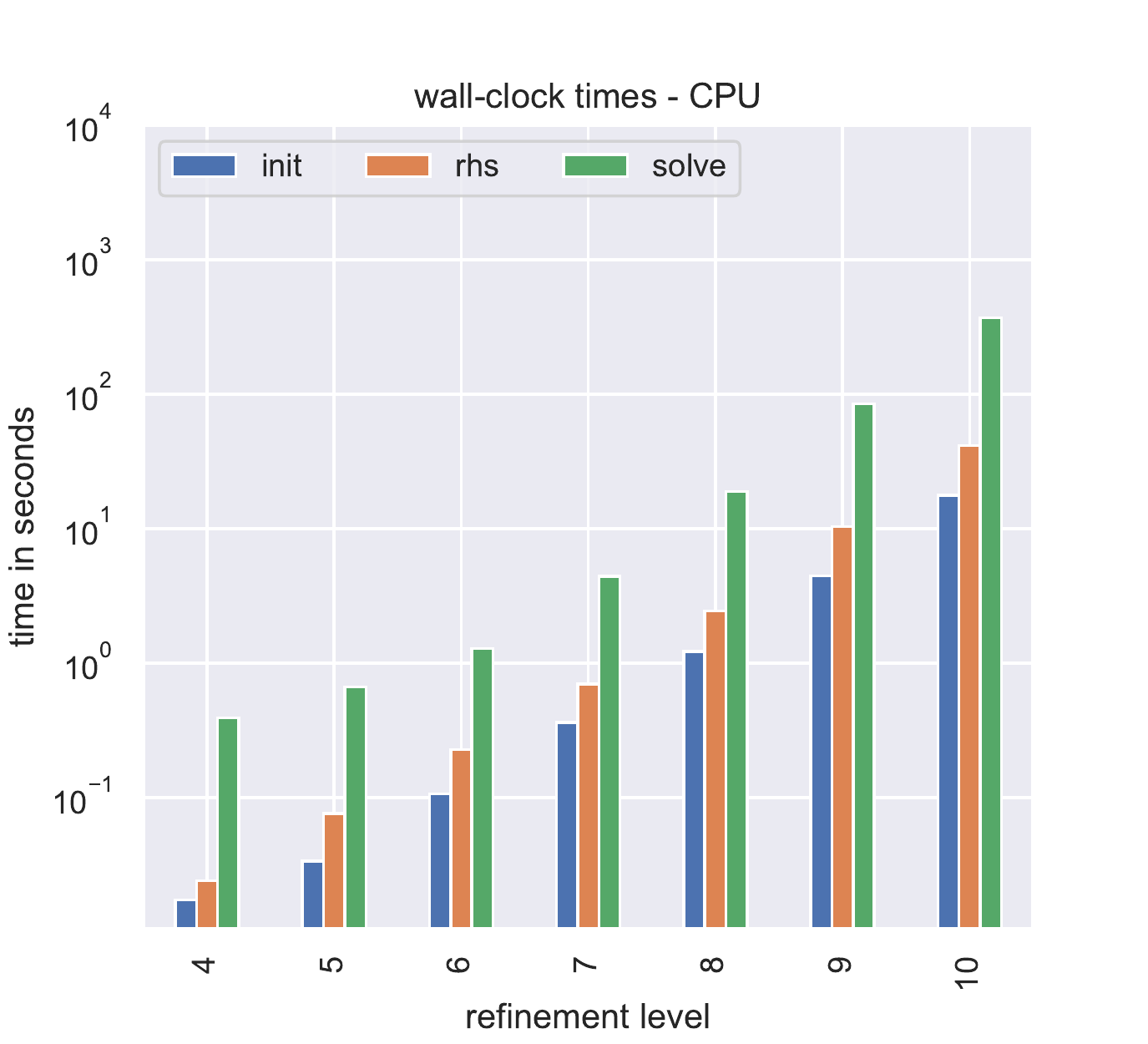}
 \includegraphics[width=0.48\linewidth]{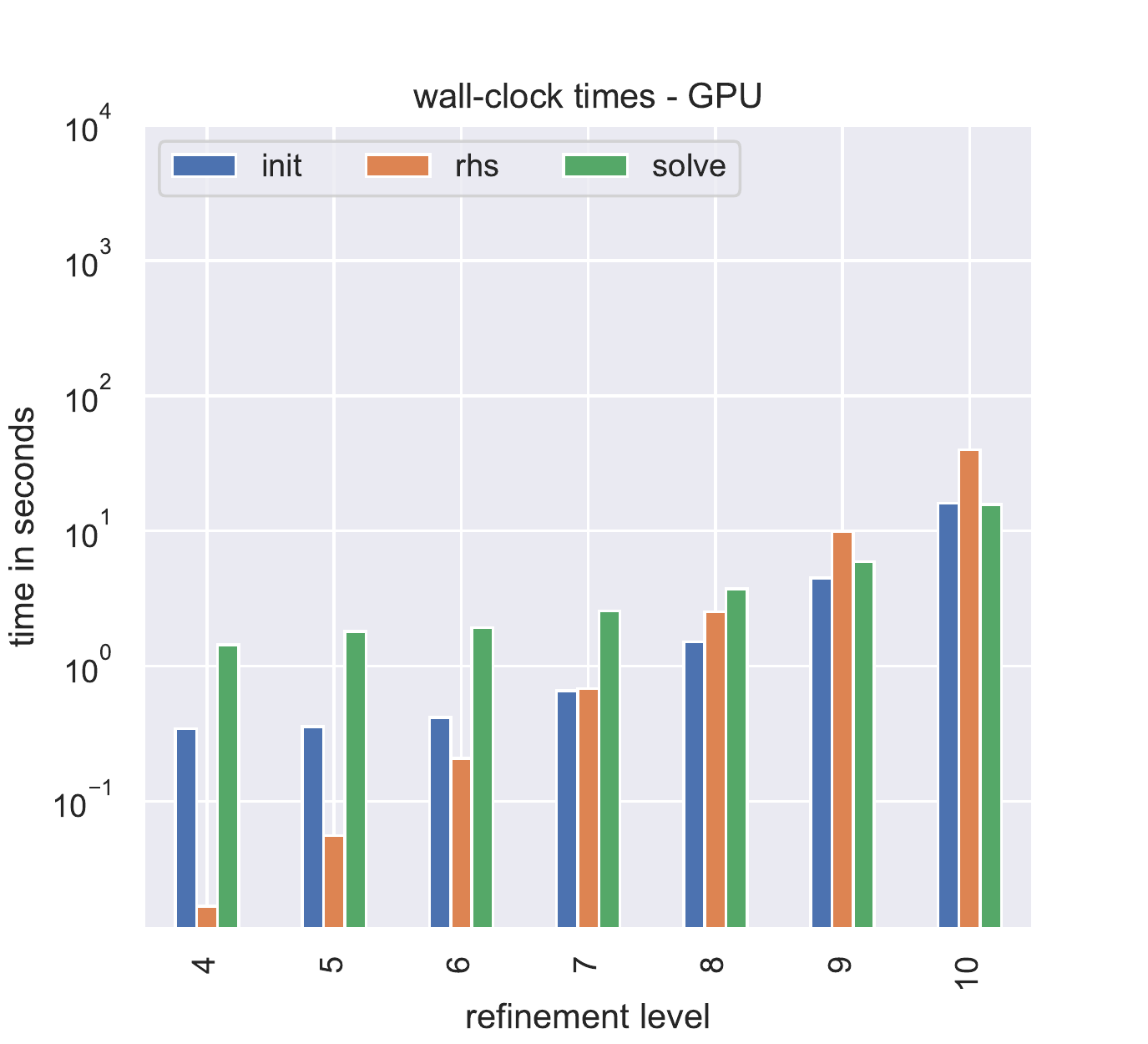}
\caption{Transport-diffusion equation: Wall-clock times CPU vs. GPU on a sequence of uniformly refined meshes.}
\label{fig:conv-dif-2}
\end{figure}

Let $\Omega \in (0,1)^2$ and $I=(0,T]$ with $T=2$. On $I\times \Omega$ we find $\theta$ subject to the following transport-diffusion problem
\begin{align}
  \partial_t\theta -\lambda\Delta\theta + (\bb \cdot \nabla)\theta \;= & \; f  \qquad\mbox{in }  (0,T]\times\Omega,\label{eq:condif_a}\\
  \theta \;= & \; \theta_b  \qquad\mbox{on }  (0,T]\times \partial\Omega \label{eq:condif_b},\\
  \theta(0,x,y) \;= & \; \theta_0 \qquad\mbox{in } 0 \times \Omega  \label{eq:condif_c}.
\end{align}
where $\lambda = 0.01$ and $\bb = (0,-1)^T$. 
The source term $f$, boundary condition $\theta_b$ and initial condition $\theta_0$ are chosen such that equations \eqref{eq:condif_a} - \eqref{eq:condif_c} have the exact solution
\[
\theta_{ex}(t,x,y) = \exp\Big(-\frac1{4} (m(t,x)^2 + m(t,y)^2) \Big),\quad
m(t,z) = \frac1{2} + \frac1{4} \cos\Big(\frac{\pi}{2} t\Big) - z .
\]
We consider the backward Euler method with a time-step size $\Delta t=0.02$.
Uniform quadrilateral spatial meshes with a range of sizes given by $h=2^{-n}$ where $n\in \{7,\dots,10\}$ are used to demonstrate the scaling efficiency of the GPU accelerated version. Hence, on the coarsest mesh with $h^{-2}=2^{2n}$ elements, the problem has $(2^7+1)^2 = 16,641$ degrees of freedom whereas the finest mesh comprises $(2^{10}+1)^2\approx 1,000,000$ degrees of freedom.
In each of the $T/\Delta t = 100$ time steps, a linear system must be solved. To ensure the robustness of the method, we use the GMRES solver, preconditioned with multigrid, to approximate the linear problems, cf. Sec.~\ref{sec:gmres}. Fig.~\ref{fig:conv-dif-1} shows the wall-clock times comparing the CPU and GPU implementations on the finest mesh level. The labels \texttt{init, rhs} and \texttt{solve} refer to the initialization (mainly assembly of system matrix), to the computation of the right hand side, and to the actual solution, respectively. Matrix assembly and the computation of the right hand side involve numerical quadrature over the mesh elements which is not ported to the GPU. Hence, no speedup is observed, GPU timings can even be higher, as copying the matrix to the GPU involves some overhead, see also Fig.~\ref{fig:conv-dif-2}. The time for solving the linear problems, however, is reduced from $\unit[304.75]{s}$ on the CPU to $\unit[9.0]{s}$ on the GPU, a factor of about 34; see also Table~\ref{tab:td} for the raw data.
The geometric multigrid solver is robust on all meshes such that the number of linear steps varies only slightly between the mesh levels.
For large problems, the non-accelerated assembly of the right hand side \texttt{rhs} hence becomes dominant. 
In Section~\ref{sec:ns} we describe how this term can be efficiently transferred to the GPU by using mass lumping.

\begin{table}[h]
\begin{center}
\resizebox{\linewidth}{!}{
\begin{tabular}{rr|rrrr|rrrr}
\toprule
    mesh&DOFs&\multicolumn{4}{c|}{CPU}    &\multicolumn{4}{c}{GPU}\\
    level & &\texttt{init}  & \texttt{rhs} & \texttt{solve}  &\texttt{sum}&\texttt{init}  & \texttt{rhs} & \texttt{solve}&\texttt{sum}\\
   \midrule
%4 & \numprint{289} &  0.0 &  0.01 &  0.41 &  0.42 &  0.27 &  0.02 &  1.17 &  1.46 \\
%5 & \numprint{1089} &  0.01 &  0.04 &  0.8 &  0.84 &  0.29 &  0.05 &  1.42 &  1.76 \\
%6 & \numprint{4225} &  0.02 &  0.14 &  1.71 &  1.87 &  0.3 &  0.17 &  1.51 &  1.98 \\
7 & \numprint{16641} &  0.07 &  0.46 &  5.33 &  5.87 &  0.35 &  0.49 &  2.05 &  2.88 \\
8 & \numprint{66049} &  0.28 &  1.66 &  20.01 &  21.96 &  0.54 &  1.69 &  2.85 &  5.08 \\
9 & \numprint{263169} &  1.29 &  7.41 &  80.63 &  89.32 &  1.41 &  7.2 &  4.47 &  13.08 \\
10 & \numprint{1050625} &  6.36 &  32.9 &  304.75 &  344.0 &  5.21 &  31.16 &  9.0 &  45.37 \\
\bottomrule
\end{tabular}}
\end{center}
\caption{Transport-diffusion problem: wall clock times on CPU and GPU in seconds. }
\label{tab:td}
\end{table}

In Figure~\ref{fig:conv-dif-2} we visually analyze the performance of the CPU (using 8 threads) and GPU version of Gascoigne 3d on the sequence of uniformly refined meshes ranging from $h=2^{-4}$ to $h=2^{-10}$. 
While the scaling of all components \texttt{init, rhs} and \texttt{solve} is linear on the CPU, the GPU-accelerated linear solver \texttt{solve} benefits from larger problems.
This is in line with theoretical considerations since the very large degree of parallelism and the very deep pipelining on the GPU requires large problems to fully utilize the computational units.
% while \texttt{init} and \texttt{rhs} are always run on the CPU.
% On very fine meshes, the computation of \texttt{rhs} is therefore dominating the overall time-stepping costs. Table~\ref{tab:td} collects all timings.

%{\color{red} I am not convinced that we should have the following paragraph. At least any reasonable reviewer %would ask that we implement it then :) Maybe refer to NS below where also RHS is assembled on the GPu.}
%We note that~\eqref{td:1}, which gets dominant for the GPU version, could easily be computed on the GPU using %mass-lumping, e.g. by approximating the current, exact computations
%\begin{equation}\label{td:1}
%b_i^f:=\big(\frac{1}{\Delta t}\theta_{n-1}+f(t_n),\phi_i\big),\quad i=1,\dots,N,
%\end{equation}
%with
%\[
%b^f \approx M_l \big( \frac{1}{\Delta t}\boldsymbol{\theta}_{n-1} + \boldsymbol{f}_n\big),
%\]
%where $\boldsymbol{\theta}_{n-1,k}=\theta_h(t_{n-1},x_k)$ and $\boldsymbol{f}_{n,k}=f(t_n,x_k)$ are the %coefficient vectors of old solution and right hand side. This simple step of optimization would drastically %reduce the computational cost for \texttt{rhs} with {\color{red} only negligible impact on the accuracy}.

%%%%%%%%

\subsection{Linear elasticity}\label{sec:elasticitynum}

\begin{figure}[t]
\begin{center}
\includegraphics[width=0.32\textwidth]{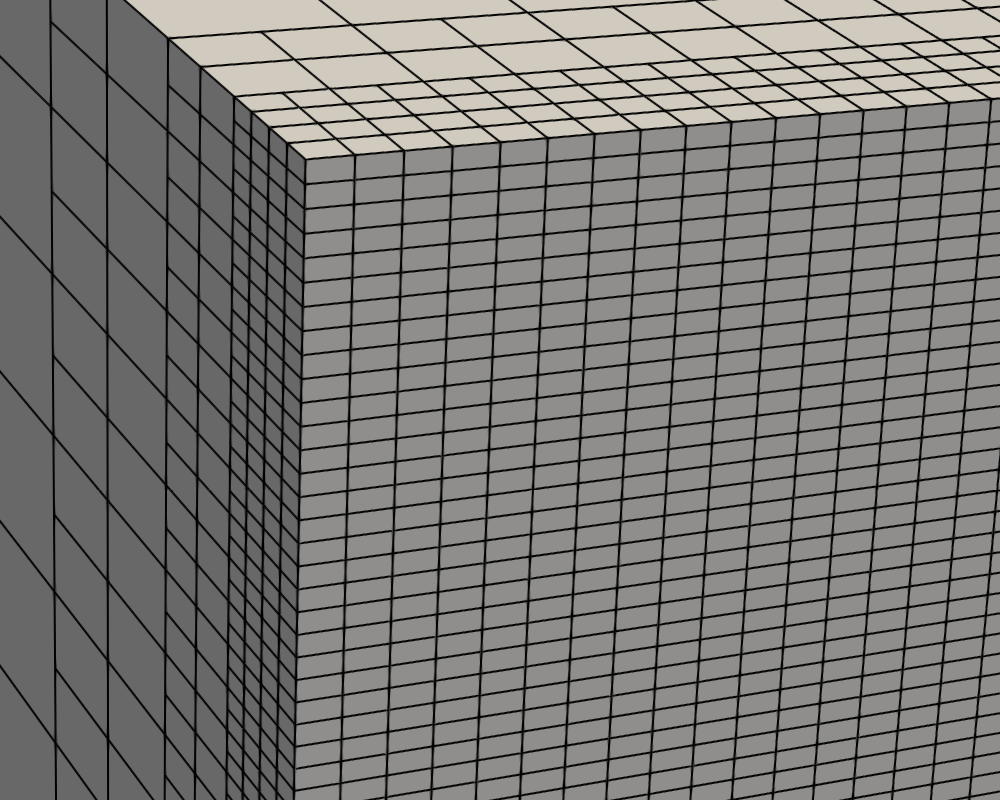}
\includegraphics[width=0.32\textwidth]{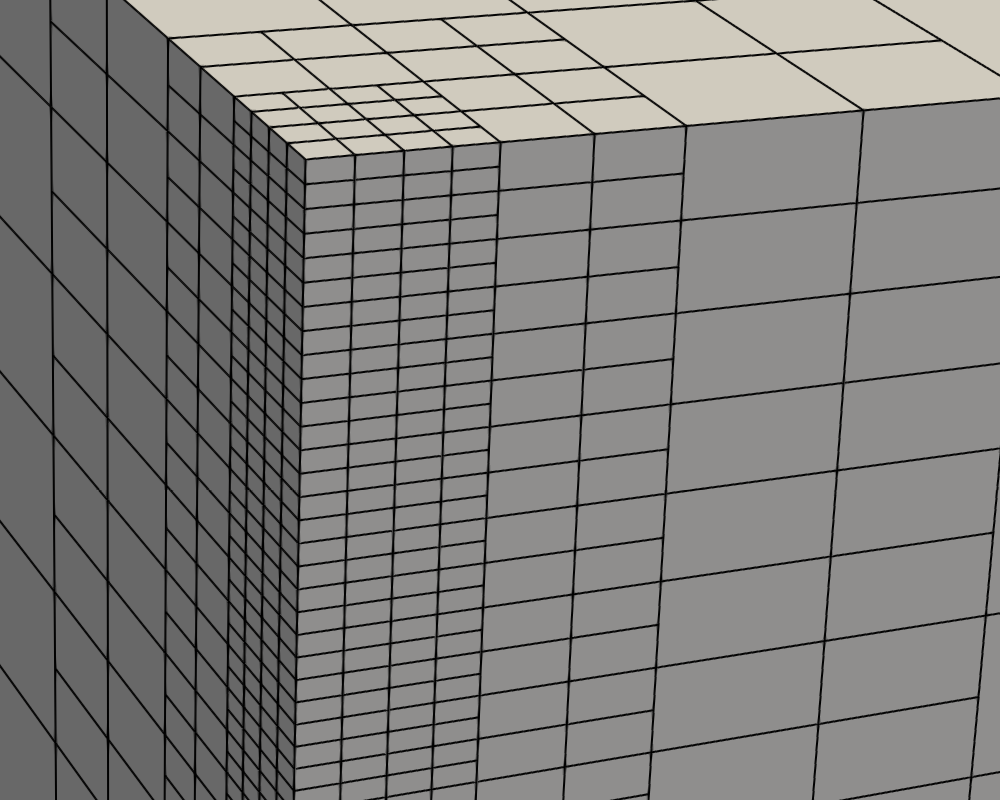}
\includegraphics[width=0.32\textwidth]{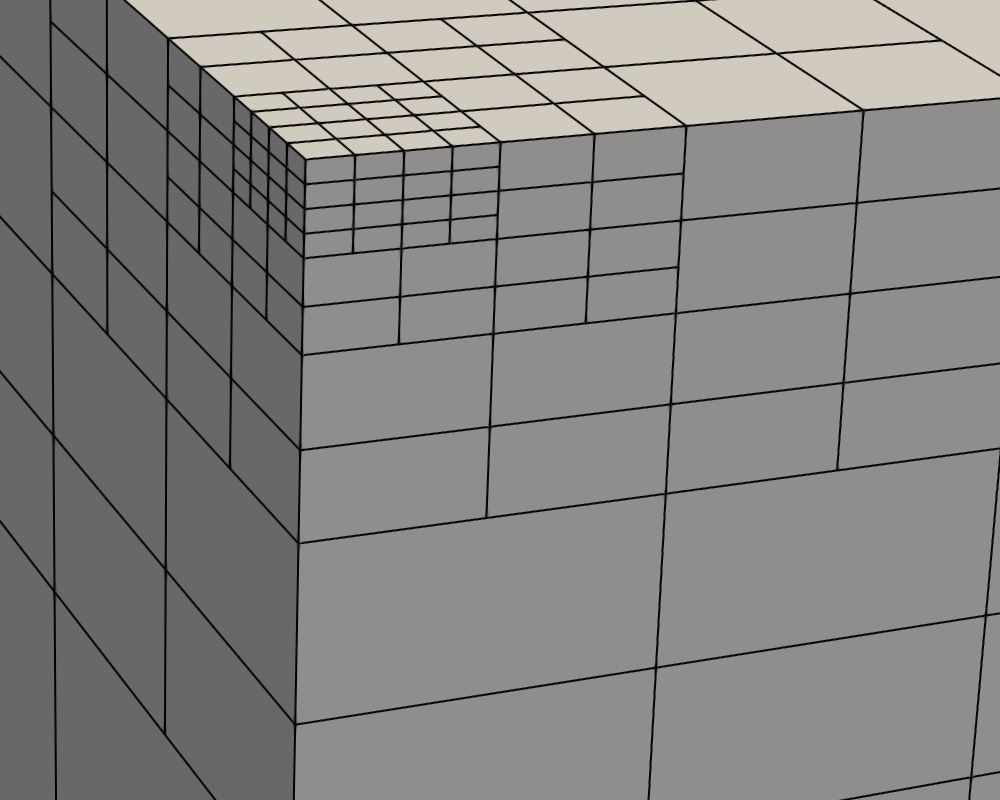}
\end{center}
\caption{Visualization of the adaptive meshes. From left to right: Refinement towards one face, towards an edge, and towards one corner. The fraction of hanging nodes serving as measure for the  unstructuredness increases from the left to the right.}
\label{fig:meshes}
\end{figure}

Second,  we consider the linear elasticity equation
\begin{align}
  \partial^2_{t}\bu - \div \bsigma(\bu) \;= \; \bf  \qquad\mbox{in }  (0,T]\times\Omega,\label{eq:elasticity}
\end{align}
on the domain $\Omega=(0,1)^3$ in the time interval $I=[0,2.5]$. By $\bu$ we denote the displacement, $\bf=(0,-1,0)^T$ is the right hand side vector, the stress tensor is
\[
\bsigma(\bu) = \lambda tr(\beps(\bu))  I + 2 \mu \beps(\bu), \quad \beps(u) = \frac1{2} (\nabla \bu + \nabla \bu^T)
\]
where $I$ is the identity tensor, $tr$ is the trace operator on a tensor and $\lambda= 8\cdot 10^4$ and $\mu = 2\cdot 10^4$ are the Lam\'e parameters. Eq. \eqref{eq:elasticity} can be written as two coupled PDEs with first order time derivatives,
\begin{align}
  \partial_{t}\bu - \bv  & \; = \; \mathbf{0}  \qquad\mbox{in }  (0,T]\times\Omega \label{eq:elasticity_1stord_a} \\
  \partial_{t}\bv - \div \bsigma(\bu) & \; = \; \bf  \qquad\mbox{in }  (0,T]\times\Omega,\label{eq:elasticity_1stord_b}
\end{align}
with homogenous Dirichlet boundary conditions and zero initial conditions for $\bu$ and $\bv$. 
We use this test case to illustrate the performance of the CUDA multigrid solver on adaptively refined meshes, see Fig.~\ref{fig:meshes} for a visualization of the used mesh. Adaptive refinement is not driven in a problem specific way here but we consider typical cases of adaptive meshing, namely refinement towards a complete face of the box, refinement towards an edge, and refinement towards a vertex. This mimics resolving singularities that have a 2d pattern (face), a 1d pattern (edge) and a 0d pattern (vertex). Table~\ref{table:ndofs} lists the number of mesh nodes and the fraction of mesh nodes that are hanging. The number of degrees of freedom is six times this number of mesh nodes as we have 3 deformation and 3 velocity components. Hanging nodes do not improve the approximation property but instead disrupt the structure of the problem by distorting the sparsity pattern of the matrix. 

\begin{table}[t]
\centering
\begin{tabular}{l | cccccc}
\toprule
Mesh levels & 1 & 2 & 3 & 4 & 5 & 6   \\
\midrule
Face&\numprint{729}	& \numprint{2925} 	& \numprint{11281}	& \numprint{43861} 	& \numprint{172505} 	& \numprint{683741} 	\\
&$0.00\%$& $4.92\%$& $6.10\%$& $6.33\%$& $6.45\%$& $6.46\%$\\
\midrule
Edge&\numprint{729} 	& \numprint{1881} 	& \numprint{4129} 	& \numprint{8569} 	& \numprint{17393} 	& \numprint{34985} 	\\
&$0.00\%$& $7.66\%$& $10.3\%$& $11.4\%$& $11.9\%$& $12.2\%$\\
\midrule
Vertex&\numprint{729} 	& \numprint{1333} 	& \numprint{1937} 	& \numprint{2541} 	& \numprint{3145} 	& \numprint{3749} 	\\
&
$0.00\%$& $8.10\%$& $11.2\%$& $12.8\%$& $13.7\%$& $14.4\%$\\
\bottomrule
\end{tabular}
\caption{Elasticity problem: Number of mesh nodes on different mesh levels (first sub-row) and fraction of nodes that are hanging nodes (second sub-row). From top to bottom we show the different adaptive refinement types, see Fig.~\ref{fig:meshes}.
\label{table:ndofs}
}
\end{table}

For time discretization of~\eqref{eq:elasticity_1stord_a}-\eqref{eq:elasticity_1stord_b} the backward Euler method is used with a step size of $\Delta t = 0.025$ resulting in $2.5/\Delta t = 100$ time steps. 
%
%% Evtl. einfach machen, wenn die Gutachter danach fragen:
%
%We use the Euler method in our experiments to simplify the numerical tests. If we were to choose the Crank-Nicolson method, for example, also a right-hand side would have to be assembled, which corresponds to matrix-vector multiplications that could be ported to the GPU. It can therefore be assumed that the GPU implementation has even greater advantages with the Crank-Nicolson method. 

\begin{figure}[t]
\begin{subfigure}{\linewidth}
\centering
  \includegraphics[width=.37\linewidth]{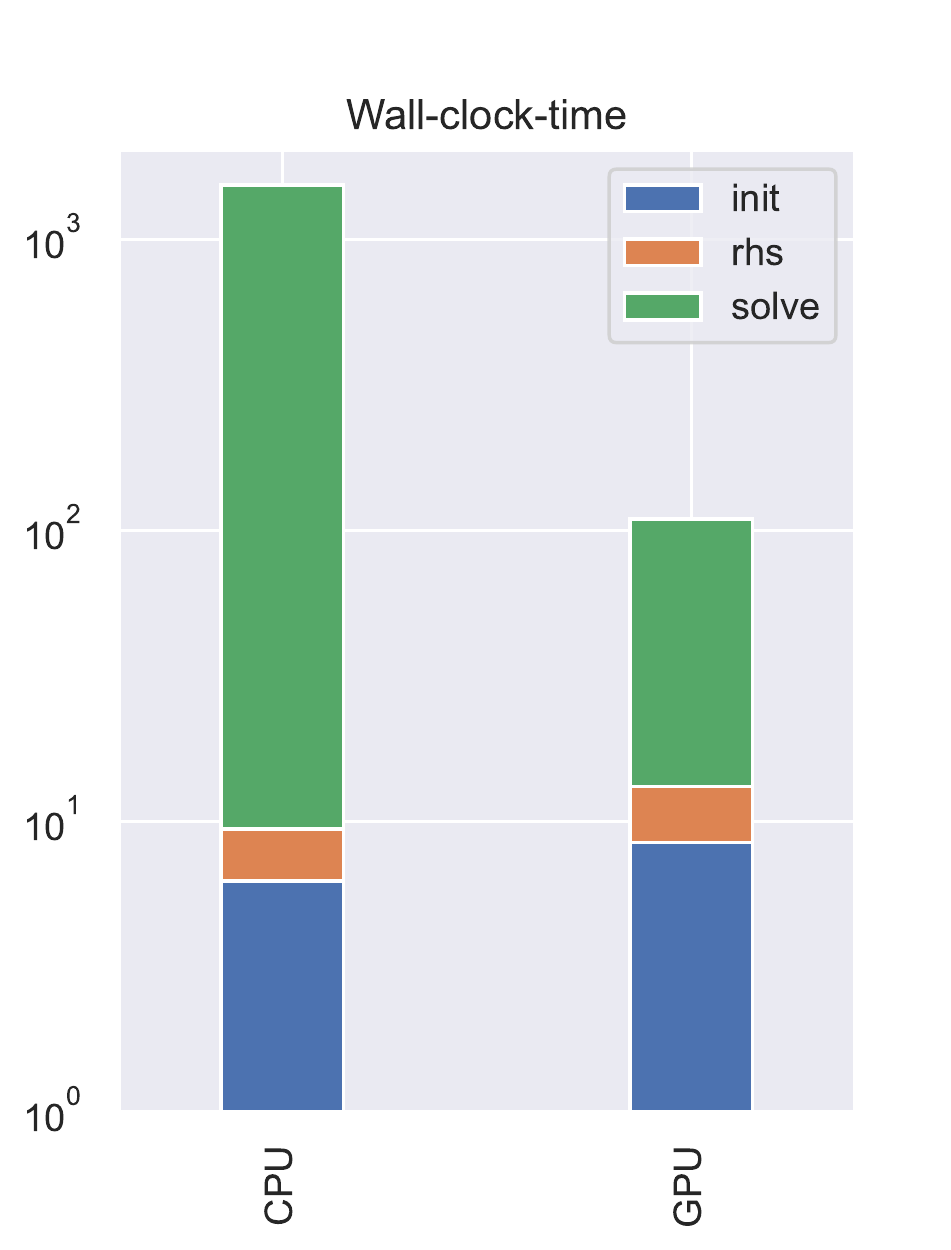}
  \includegraphics[width=.37\linewidth]{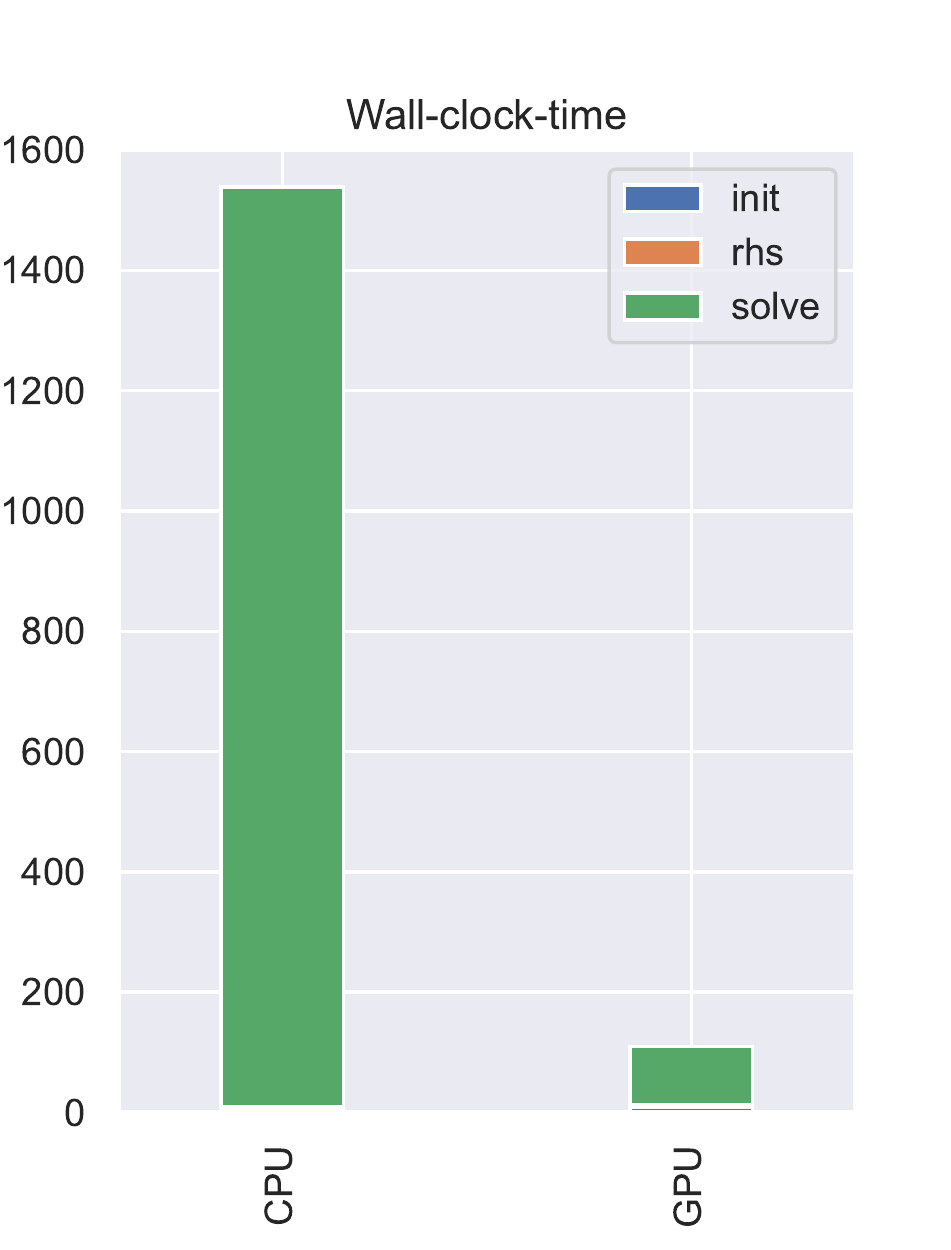}
  \subcaption{Refinement towards edge on mesh level $6$, logarithmic scale (left) and linear scale (right).}
  \label{fig:elasticity-1}
  \end{subfigure}
  \begin{subfigure}{\linewidth}
    \centering
   \includegraphics[width=.48\linewidth]{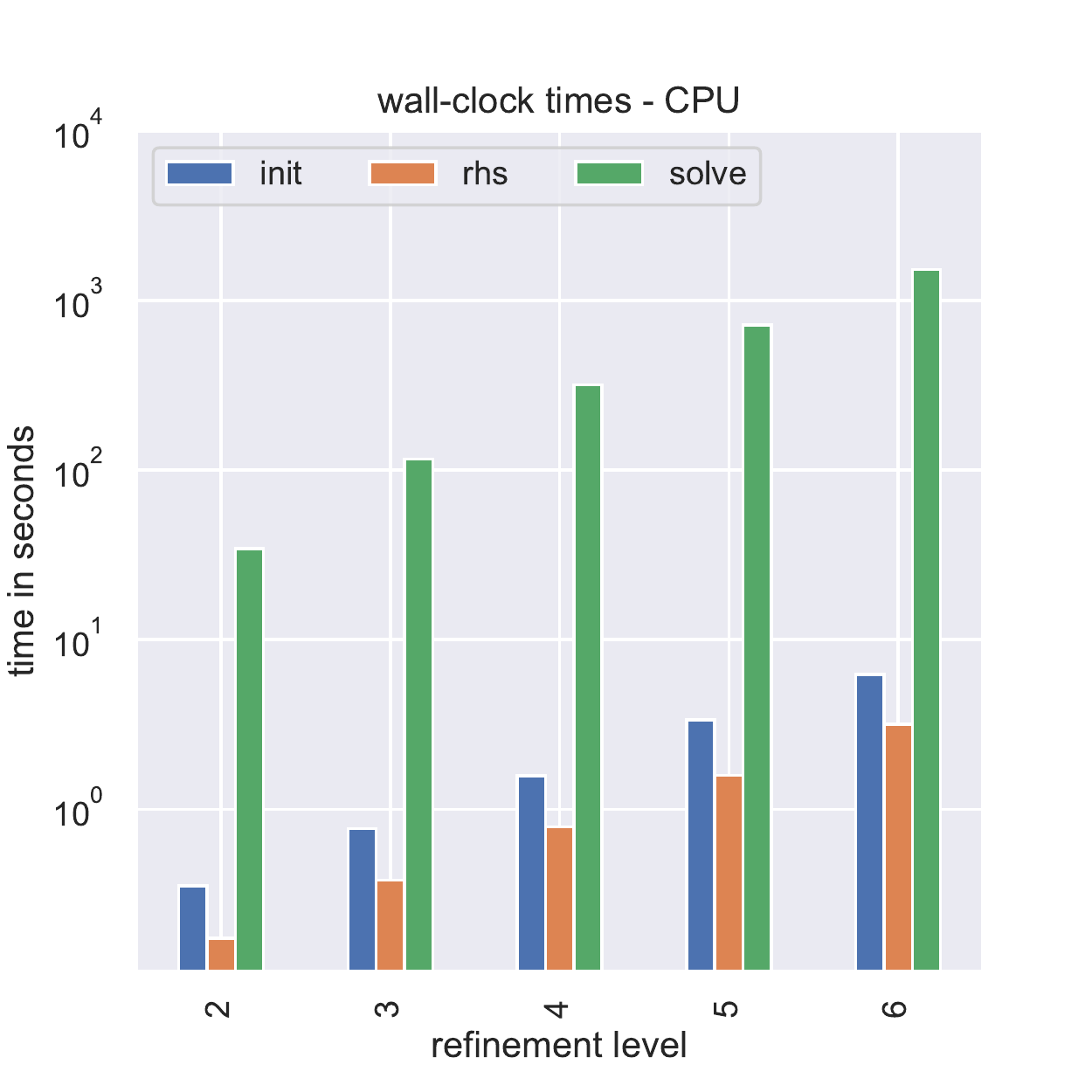}
   \includegraphics[width=.48\linewidth]{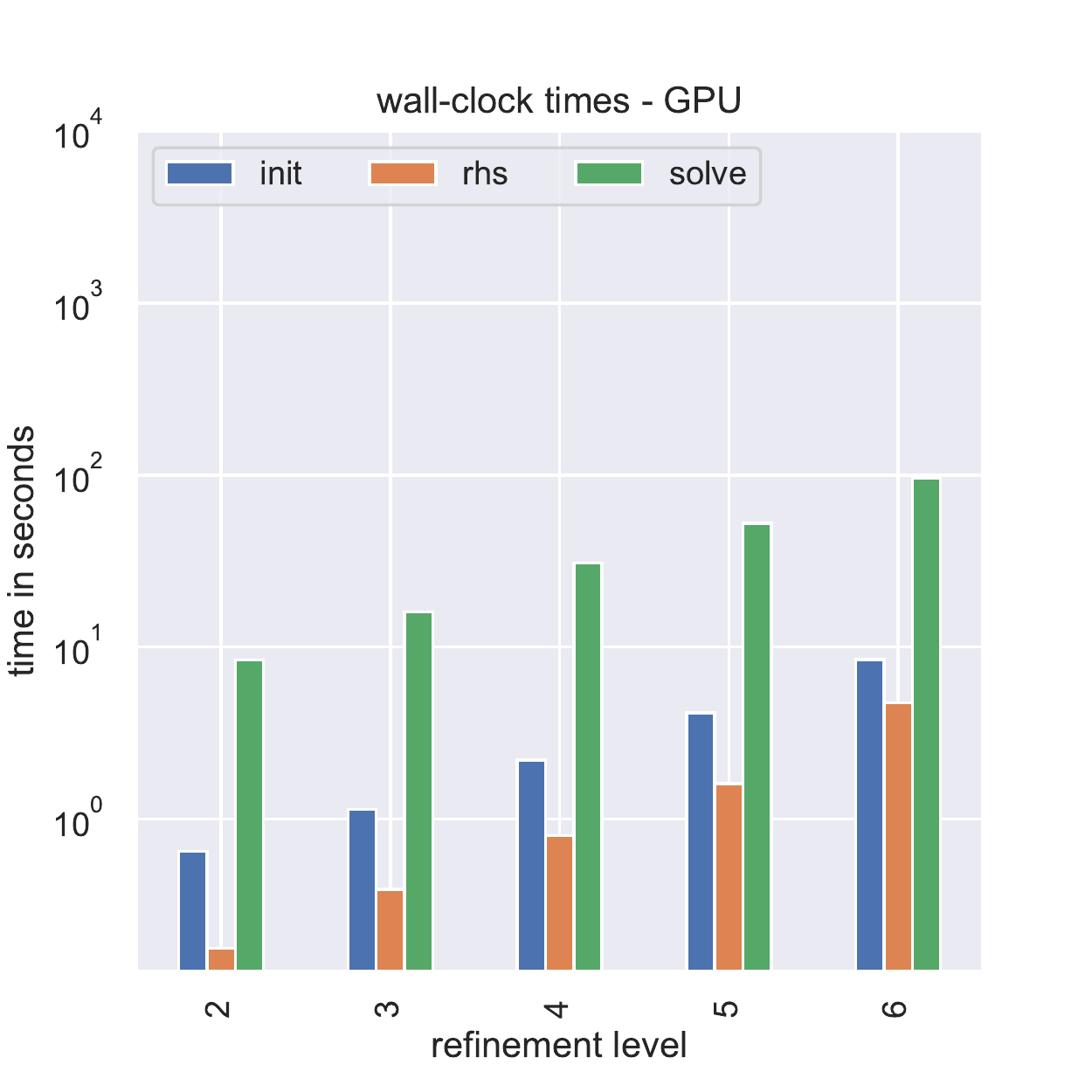}
   \caption{Refinement towards edge on a sequence of locally refined meshes.}
  \label{fig:elasticity-2}
    \end{subfigure}

    \caption{Elasticity problem. Wall-clock times on CPU and GPU.}
\end{figure}

We provide a detailed discussion of the intermediate case, the refinement towards one edge. Fig.~\ref{fig:elasticity-1} gives the timings on the finest mesh level (34\,985 nodes, hence 209\,910 degrees of freedom, about 25\,600 of them in hanging nodes) for the CPU version and the GPU version. As for the transport-diffusion equation, times for initialization \texttt{init} and right hand side \texttt{rhs} do not change, as these parts are not implemented on the GPU.  
The computation time for the linear solver is, however, drastically reduced, see  Table~\ref{tab:elasticity:times:edge}. 
On the finest mesh, the solution time reduces from $\unit[306]{s}$ to $\unit[15]{s}$, i.e. by a factor of 20.
As the time for matrix assembly (\texttt{init}) and computing the right hand side (\texttt{rhs}) are negligible, the overall speedup is still close to \blue{17} (compared to about \blue{8} for the transport-diffusion problem on a much finer mesh).  As for the transport-diffusion problem, the initialization (\texttt{init}) may require more time on the GPU as this includes the cost for copying the matrix from the CPU to the GPU. Also, since the effort is negligible for this problem, the right hand side has been assembled on the CPU for simplicity. Therefore, the timinigs on the GPU include overhead for copying the data to the GPU.

Fig.~\ref{fig:elasticity-2} shows the scaling of the implementation with respect to the refinement. Note that the meshes are refined locally. Therefore, the number of multigrid levels increases in each step but the overall number of unknowns grows only slowly (not by a factor of 8 that would be expected on a 3d hexahedral mesh). Hence more and more unstructured features and hanging nodes appear on higher mesh levels. Due to the slow growth in the number of degrees of freedoms, the run-times on CPU and GPU increase only slowly with each refinement level, see also Table~\ref{table:ndofs}. The results show very good efficiency of the GPU multigrid solver on locally refined meshes with no negative impact of having more and more hanging nodes and less regular structure in the problem.

Next we compare the performance of the GPU implementation on the different types of adaptive mesh shown in Fig.~\ref{fig:meshes} with number of unknowns indicated in Table~\ref{table:ndofs}. The complexity of these mesh types is very different. While the number of nodes increases like $4^l$ ($l$ is the mesh level) for face refinement, which is typical for 2d problems, and $2^l$ for the edge refinement, the typical behavior for 1d problems, the scaling is much slower for the vertex-refinement.

%\FloatBarrier

\begin{figure}[t]
      \centering
  \includegraphics[width=.32\linewidth]{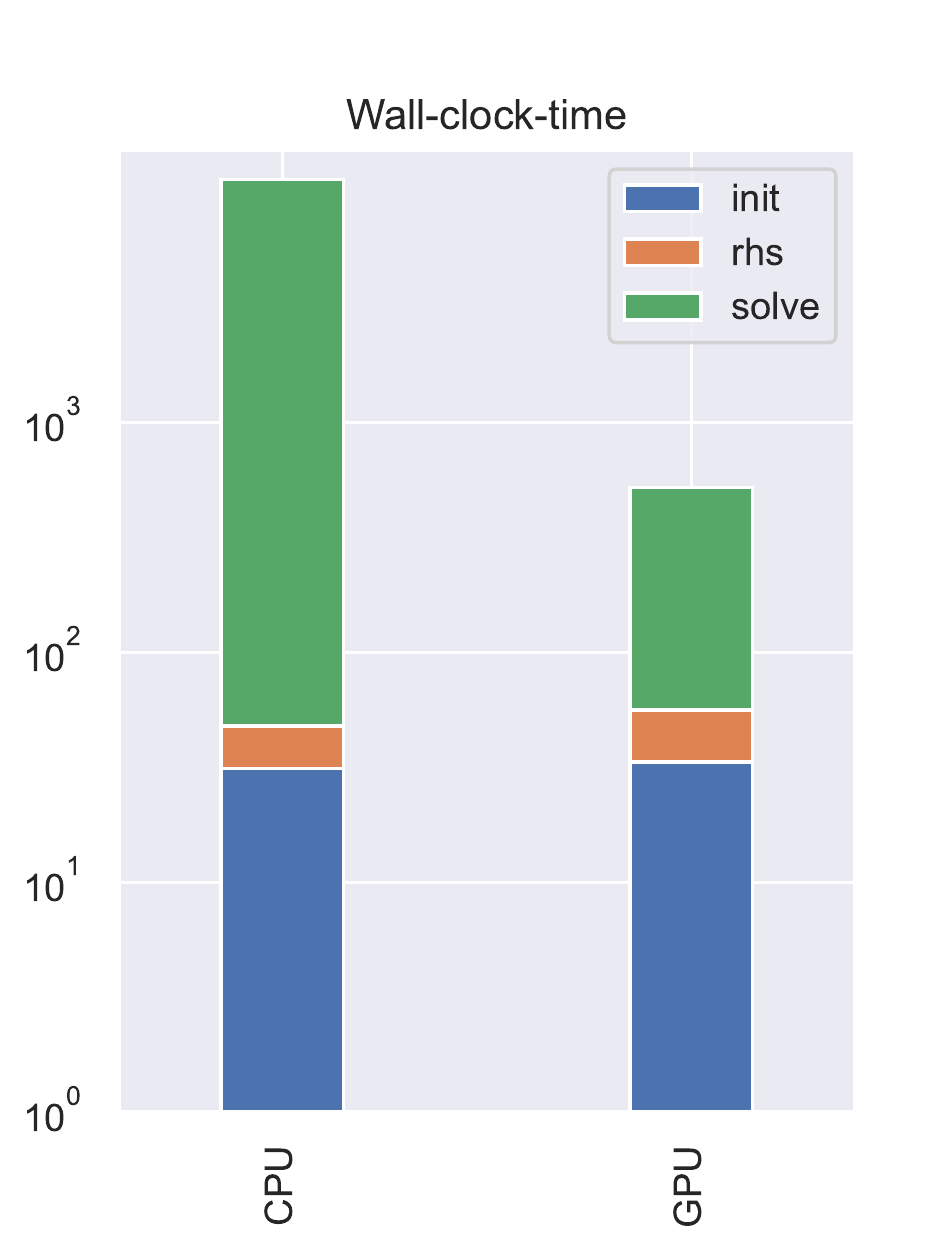}
  \includegraphics[width=.32\linewidth]{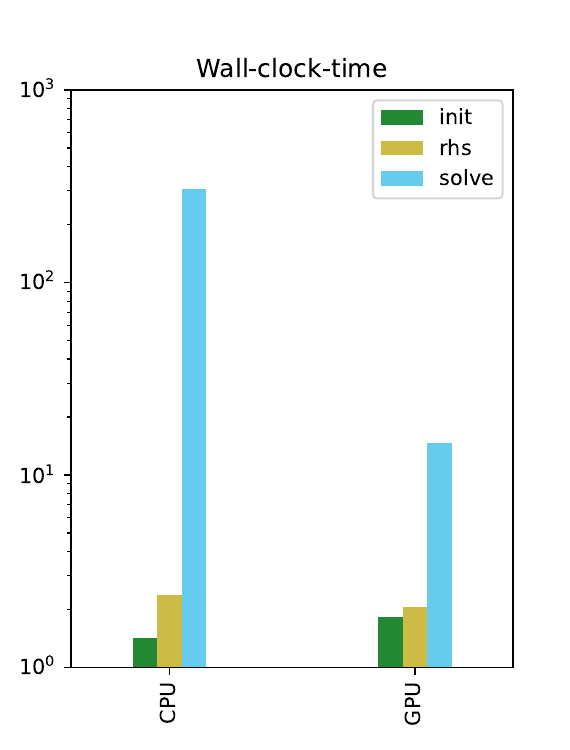}
  \includegraphics[width=.32\linewidth]{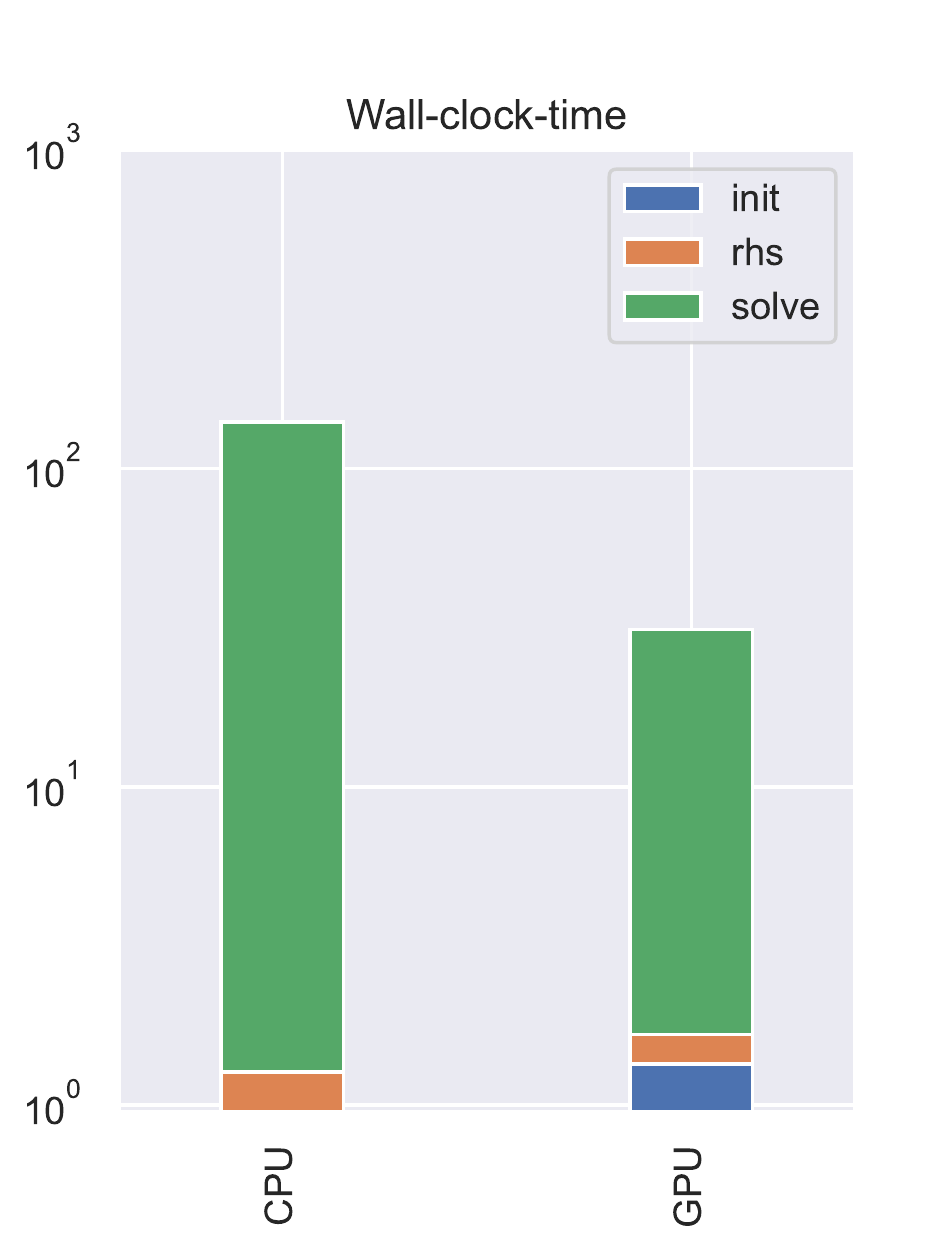}
  \caption{Comparison of the refinement types (from left to right): 5 refinements towards one face (left), 6 refinements towards one edge (middle), 6 refinements towards one vertex (right).}
\label{fig:elasticity-3}
\end{figure}

\begin{table}[t]
\begin{subtable}{\linewidth}
\centering
\resizebox{\linewidth}{!}{
\begin{tabular}{rr| rrrr|rrrr}
\toprule
mesh & DOFs&\multicolumn{4}{c|}{CPU}& \multicolumn{4}{c}{GPU}\\
level  && \texttt{init} & \texttt{rhs} & \texttt{solve} &\texttt{sum}  & \texttt{init}  & \texttt{rhs} & \texttt{solve} &\texttt{sum}\\
 \midrule
2 & \numprint{17550} &  0.17 &  0.22 &  41.78 &  42.17 &  0.45 &  0.34 &  7.52 &  8.32 \\
3 & \numprint{67686} &  0.7 &  0.89 &  174.96 &  176.55 &  0.99 &  1.27 &  14.65 &  16.92 \\
4 & \numprint{263166} &  2.99 &  3.56 &  1031.13 &  1037.68 &  3.45 &  2.98 &  38.22 &  44.64 \\
5 & \numprint{1035030} &  11.91 &  14.01 &  4431.94 &  4457.85 &  12.56 &  10.78 &  122.57 &  145.9 \\
\bottomrule
\end{tabular}}
\subcaption{Refinements towards a face of the domain.}
\label{tab:elasticity:times:face}
\end{subtable}
\begin{subtable}{\linewidth}
\centering
\resizebox{\linewidth}{!}{
\begin{tabular}{rr| rrrr|rrrr}
\toprule
mesh & DOFs&\multicolumn{4}{c|}{CPU}& \multicolumn{4}{c}{GPU}\\
level  & &\texttt{init} & \texttt{rhs} & \texttt{solve} &\texttt{sum}  & \texttt{init}  & \texttt{rhs} & \texttt{solve} &\texttt{sum}\\
 \midrule
%2 & \numprint{11286} &  0.06 &  0.14 &  14.43 &  14.62 &  0.34 &  0.22 &  3.97 &  4.52 \\
3 & \numprint{24774} &  0.13 &  0.29 &  38.22 &  38.64 &  0.43 &  0.49 &  5.93 &  6.85 \\
4 & \numprint{51414} &  0.3 &  0.64 &  82.92 &  83.86 &  0.61 &  0.99 &  7.64 &  9.24 \\
5 & \numprint{104358} &  0.66 &  1.19 &  165.29 &  167.14 &  1.0 &  1.62 &  10.38 &  13.01 \\
6 & \numprint{209910} &  1.42 &  2.38 &  305.58 &  309.39 &  1.82 &  2.07 &  14.7 &  18.59 \\
\bottomrule
\end{tabular}}
\subcaption{Refinement towards an edge of the domain.}
\label{tab:elasticity:times:edge}
\end{subtable}
\begin{subtable}{\linewidth}
\centering
\resizebox{\linewidth}{!}{
\begin{tabular}{rr| rrrr|rrrr}
\toprule
mesh&DOFs & \multicolumn{4}{c|}{CPU}& \multicolumn{4}{c}{GPU}\\
level&  & \texttt{init} & \texttt{rhs} & \texttt{solve} &\texttt{sum}  & \texttt{init}  & \texttt{rhs} & \texttt{solve} &\texttt{sum}\\
 \midrule
%2 & \numprint{7998} &  0.04 &  0.09 &  8.41 &  8.54 &  0.32 &  0.16 &  3.58 &  4.06 \\
3 & \numprint{11622} &  0.07 &  0.13 &  13.04 &  13.24 &  0.37 &  0.22 &  4.73 &  5.33 \\
4 & \numprint{15246} &  0.11 &  0.17 &  17.68 &  17.96 &  0.42 &  0.31 &  6.12 &  6.85 \\
5 & \numprint{18870} &  0.16 &  0.21 &  22.7 &  23.07 &  0.48 &  0.36 &  7.09 &  7.94 \\
6 & \numprint{22494} &  0.21 &  0.25 &  28.21 &  28.67 &  0.56 &  0.42 &  7.9 &  8.88 \\
\bottomrule
\end{tabular}}
\subcaption{Refinement towards a vertex.}
\label{tab:elasticity:times:vertex}
\end{subtable}
\caption{Elasticity problem: measured values for the wall clock times on CPU and GPU on different refinement types.}
\end{table}

Figure~\ref{fig:elasticity-3} shows that the efficiency of the GPU parallelization is robust. On the finest mesh level belonging to the vertex-case in Table~\ref{tab:elasticity:times:vertex}, the computational time of the solver only drops from $\unit[29]{s}$ on the CPU to $\unit[8]{s}$ on the GPU, a factor of about 3.5. However, this problem has just \numprint{25000} unknowns and is too small to properly utilize a GPU. It is important to note that the GPU parallelization does not lead to an overhead at any point, which negatively affects the entire computing time. 

\subsection{Navier-Stokes equations}\label{sec:ns}

As third example, we consider the time-dependent Navier-Stokes equations in their dimensionless form
\begin{align}
    \label{eq:nstokes1}
  \partial_t\bu -\frac1{\Re}\Delta\bu+ \div  (\bu \otimes \bu ) + \nabla p \; &= \;0  \qquad\mbox{in }  (0,T]\times\Omega,
  \\
    \div\bu \; &= \; 0 \qquad\mbox{in }(0,T]\times\Omega .
    \label{eq:nstokes2} 
\end{align}
By $\bu: (0,T] \times \Omega \to \mathbb{R}^3$ we denote the velocity field and by $p: (0,T] \times \Omega \to \mathbb{R}$ the scalar pressure.
In our numerical experiment, we use the driven cavity benchmark problem on the time interval $I=(0,12]$ and on the three dimensional domain $\Omega=(0,1)\times (0,1)\times (0,2)$. 
We choose the Reynolds number $\Re = 10^3$ and the boundary conditions
$$
\bu(t,x,y,z) = 
\begin{cases}
  (0,1,0)^T & \text{if } x=1,\\
  (0,0,0)^T & \text{otherwise.}
\end{cases}
$$
Homogenous initial condition $\bu(0,x,y,z) = \bu_0 := \mathbf{0}$ hold at time $t=0$. 
We set $\bV = [H^1_0(\Omega)]^d$ and $Q:=\{q \in L^2(\Omega)\,:\,\int_\Omega q\,\text{d}x = 0\}$. The weak formulation of the time-dependent Navier-Stokes equations in its semi-discrete form is 
\begin{alignat}{2}
  (\partial_t\bu, \bchi) +\frac1{\Re}(\nabla\bu,\nabla\bchi) -(\bu \otimes \bu,\nabla \bchi) -  (p,\div \bchi) \;= \;& (\bf,\bchi) \qquad &\forall \bchi \in \bV,\label{eq:weaknstokes1}\\
    (\div\bu,\xi) \;= \; & 0\qquad &\forall \xi \in Q.\label{eq:weaknstokes2}
\end{alignat}
Here, we used the divergence form of the convective term once integrated, i.e.
\[
\big(\bu\cdot\nabla\bu,\bchi\big) = 
\big(\div(\bu\otimes\bu),\bchi\big) = 
-\big(\bu\otimes\bu,\nabla\bchi\big) .
\] 
This representation will be crucial for an efficient realization on the GPU, since in the discrete setting $\big(\bu\otimes\bu,\nabla\bchi\big)$ can be approximated by sparse matrix vector products, as will be explained shortly.

\begin{algorithm}[t]
  \small
  \caption{The fully discrete solution procedure for the Navier-Stokes equations}
  \label{algo:2}
\rule{\textwidth}{\algow}
Given $\bu_0$. Set $q_h^0 := 0$ and choose $p_h^0 \in Q_h$ s.t. \eqref{eq:init_pres} is fulfilled. For $m=1, \dots, N$ calculate
\begin{align*}
\intertext{\textbf{Step 1:}  
    Find $\bu_h^m \in \bV_{h}$ such that}
  \frac1{k}(\bu_h^m, \bchi)_{*} & \; = \; (\bf^{m-1}+\frac1{k}\bu^{m-1}_h,\bchi)_*
- A((\bu^{m-1}_h, p^{m-1}_h + q^{m-1}_h);(\bchi,0))
  \quad \forall \bchi \in \bV_{h}.
  \label{eq:pre_discrete}
\intertext{\textbf{Step 2:}  Find $q^m_h \in S_h$ such that}
(\nabla q^m_h, \nabla  \varphi) & \; = \; -\frac1{k}(\div \bu_{h}^m, \varphi) \quad \forall \varphi \in S_{h}.
\intertext{\textbf{Step 3:}  Find $p^m_{h} \in Q_{h}$ such that }
( p_{h}^m,\varphi)_{*}  &\;=\; (p_{h}^{m-1}+q_{h}^m,\varphi)_*-(\nu\,\div\,\bu_{h}^m,\varphi)  \quad  \forall \varphi \in Q_{h}.
\end{align*}
% \rule{\textwidth}{\algow}
\end{algorithm}

For the spatial discretization of \eqref{eq:weaknstokes1}-\eqref{eq:weaknstokes2}, we employ finite elements with the inf-sup stable $\mathbb{Q}_1$-iso-$\mathbb{Q}_2/\mathbb{Q}_1$ pair. 
We consider a shape regular mesh $\Omega_h$ of the domain and its equidistant refinement $\Omega_{h/2}$. The discrete spaces $\bV_h \subset \bV$ and $Q_h \subset Q$ consist of $\mathbb{Q}_1$ polynomials in each cell of $\Omega_{h/2}$ and $\Omega_h$, respectively. The semi-discrete problem is obtained with $\bu_h,\bchi_h \in \bV_h$ and $p_h,\xi_h \in Q_h$.  
Application of an implicit time discretization for this semi-discrete problem requires fixed point iterations due to the nonlinear transport term. Each iteration step for linearization requires thereby to resolve a nonsymmetric, indefinite system matrix owing to the transport term and the saddle point structure that arises from velocity pressure coupling. Here, we follow an semi-implicit approach where the momentum equation is discretized explicitly.

Due to better stability and the advantageous diagonal structure, we employ a lumped mass matrix, i.e. we approximate
\begin{equation}\label{mass_lumping}
(\bv_h,\phi_i)_{i=1}^N \approx M_h^l\bv,
\end{equation}    
where $\phi_i^h : \Omega \to \mathbb{R}$ are the $\mathbb{Q}_1$ basis functions, $\bv\in\mathbb{R}^{N_{nodes}\times 4}$ is the coefficient vector, and $M_h$ denotes the lumped mass matrix given by
\[
M_h^l = \operatorname{diag}(m_1,\dots,m_{N_{dof}}),\quad
m_i = \sum_j (\phi_j,\phi_i).
\] 
Moreover, we apply a mass lumping approach also in the nonlinear term, which enables us to implement the convective term as a matrix-vector product with a pre-assembled sparse matrix, see also~\cite{KR2024}. 
In particular, let $\mathcal{I}_h: \bV \to \bV_h$ be the nodal interpolation operator. For each $\bv_h,\; \bchi_h \in \bV_h$, we approximate 
\begin{align}  \label{conv_mass_lumping}
  \begin{split}
(\div (\bv_h \otimes \bv_h), \bchi_h) =& - (\bv_h \otimes \bv_h, \nabla \bchi_h) \\
\approx &  - (\mathcal{I}_h (\bv_h \otimes \bv_h), \nabla \bchi_h)  =: n(\bv_h\otimes \bv_{h}, \bchi_h).
  \end{split}    
\end{align}
In finite element notation, this discrete approximation amounts to
\[
\mathcal{I}_h (\bv_h \otimes \bv_h)
=\sum_i (\bv_i\otimes \bv_i)\phi_i^h
\approx \Big(\sum_i \bv_i\phi_i^h\Big)\otimes\Big(\sum_j \bv_j\phi_j^h\Big).
\]
Note that in the interpolation nodes it holds that
\[
\mathcal{I}_h(\bv_h\otimes \bv_h)(\bx_k)=\bv_h(\bx_k) \otimes \bv_h(\bx_k).
\]
The error caused by mass lumping in the zeroth-order term \eqref{mass_lumping} and the convective term \eqref{conv_mass_lumping} are of the same order as the polynomial approximation error and therefore do not affect the error asymptotics. This approach is sometimes denoted as `fully practical finite element method'~\cite{Barrett2001} as it allows for error estimates that reflect the full error including numerical quadrature.
We have recently published the complete error analysis for an explicit pressure correction method based on this approximation~\cite{KR2024}. Both theoretical analysis and numerical examples show that the explicit treatment based on the interpolation of the nonlinearity does not lead to any additional error. 

For the convective term we assemble three sparse matrices
\[
C_{h,\{x,y,z\}} = (\phi_j,\partial_{\{x,y,z\}}\phi_i)_{i,j=1}^{N_{nodes}},
\]
for the three components $x$, $y$, $z$ and compute the three vectors
\begin{equation}\label{tp}
\bv^1_{i,c} = \bu_{i,1}\bu_{i,c},\;
\bv^2_{i,c} = \bu_{i,2}\bu_{i,c},\;
\bv^3_{i,c} = \bu_{i,3}\bu_{i,c},\; .
\end{equation}
for $c=1,2,3,\; i=1,\dots,N_{nodes} $.
The residual can then be evaluated using matrix-vector products as
\begin{equation}\label{multC}
-\big(\bu_h\otimes \bu_j,\nabla\bchi\big)_{i=1}^{N_{nodes}}\approx
- C_{h,x} \bv^1 - C_{h,x} \bv^2 - C_{h,x} \bv^3.
\end{equation}
The computation of the component-wise products in~\ref{tp} is not a standard operation that can be expressed in linear algebra. We will therefore require custom kernels to assemble these vectors efficiently on the GPU.

\subsubsection{Explicit projection solver}
Let $k>0$ be a time step and  $N=T/k$. We set $t_n := n\cdot k$ for $0\leq n \leq N$ and 
\begin{align*}
  A((\bu, p);(\bchi,\xi)) := 
  - \nu(\nabla \bu, \nabla \bchi)
  + n(\bu\otimes \bu, \bchi) -(p,\div \bchi) + (\xi,\div \bu).
\end{align*}
The fully discrete solution algorithm we employ solves the momentum equation explicitly and updates the pressure field by solving a Poisson problem (see Alg. \ref{algo:2}). Hence, we introduce the solution space $S=H^1(\Omega)$ and its discrete counterpart $S_h$ that consists of $\mathbb{Q}_1$ polynomials on each cell $T$ of $\Omega_h$.
This type of predictor-corrector methods for approximating incompressible flows requires an initial pressure field which must be calculated from the Poisson equation if the right hand side and initial velocity are not zero
\begin{align}\label{eq:init_pres}
  (\nabla p^0_h, \nabla  \varphi) & \; = \; -(f^0, \nabla \varphi) + ( (\bu_0 \cdot \nabla) \bu_0, \nabla\varphi) +( \nu\Delta \bu_0, \varphi)_{\Gamma} \quad \forall \varphi \in S_{h}.
\end{align}
For justification and well-posedness of this equation, see \cite{rannacher1982}. Moreover, note that the velocity solutions $\bu_h^m$ are so-called predictor solutions. The corrector velocities were eliminated as described in \cite{guermond1996}. The complete error analysis for this fully explicit variant of a pressure correction scheme is found in~\cite{KR2024} and its extension to the Boussinesq approximation is described in~\cite{RKR2024}.

We set $\Delta t = 10^{-4}$ and consider $(N_x,\; N_y, \; N_z) = (32,32,64)$ elements in each coordinate direction.
The mesh coordinates $(x_i,\;y_j,\;z_k)$ for $i \in \{0,\cdots,N_x\},$  $j \in \{0,\cdots,N_y\}$ and $k \in \{0,\cdots,N_z\}$ are defined as
\[
\begin{aligned}
  x_i & =  \frac1{2} \Big ( 1 + \cos \big(\frac{i \cdot \pi}{N_x}\big) \Big), &\qquad i &\in \{0,\cdots,N_x\},\\
  y_j & =  \frac1{2} \Big ( 1 + \cos \big(\frac{j \cdot \pi}{N_y}\big) \Big),& \qquad j &\in \{0,\cdots,N_y\},\\
  z_k & =  1 + \sin \Big ( \frac{(k-N_z)\cdot \pi}{2N_z} \Big), &\qquad k &\in \{0,\cdots,N_z\}.
\end{aligned}
\]
The meshes contain anisotropic elements. Therefore, to obtain better robustness, we again embed the geometric multigrid solver as preconditioner in the GMRES method.

\subsubsection{GPU realization}\label{sec:gpu:custom}

The explicit pressure correction scheme in Algorithm~\ref{algo:2} consists largely of the same parts that have already been used for the nonstationary linear problems that have been presented above. The multigrid method is fully run on the GPU as well as all parts of the GMRES solver that work with global matrices and vectors scaling with the dimension of the finite element space. \textbf{Step 1} and \textbf{Step 3} are explicit and only require sparse matrix vector products as well as inversion of the lumped mass matrix. This however is already stored as a diagonal matrix containing the inverse elements such that standard cuSPARSE methods can be used.

An operation that is not easily expressed in cuBLAS or cuSPARSE is the local assembly of the outer products $\mathbf{u}_i\otimes\mathbf{u}_i\in\mathbb{R}^{3\times 3}$ in~\eqref{tp}. 
To understand the effect this has on the performance, the numerical results below show two versions: one, simply denoted by (GPU), where a custom CUDA kernel is used to evaluate~\eqref{tp} directly on the GPU and an intermediate version, denoted as (GPU$^*$), where we compute these products on the CPU. This requires transferring the solution vector back and forth whenever this matrix-vector product is required. We add this intermediate results to demonstrate the importance of minimizing memory transfers.

We also use custom kernels to efficiently compute the products with the rectangular matrices corresponding to the discretization of the gradient of the pressure and the divergence of the velocity that is required in preprocessing. Appendix~\ref{app:impl} gives details on the implementation and in particular on the required changes in the user code to employ GPU parallelization and on the integration of custom kernels.

\subsubsection{Numerical Results}

\begin{figure}[t]
\centering
\includegraphics[width=.49\linewidth]{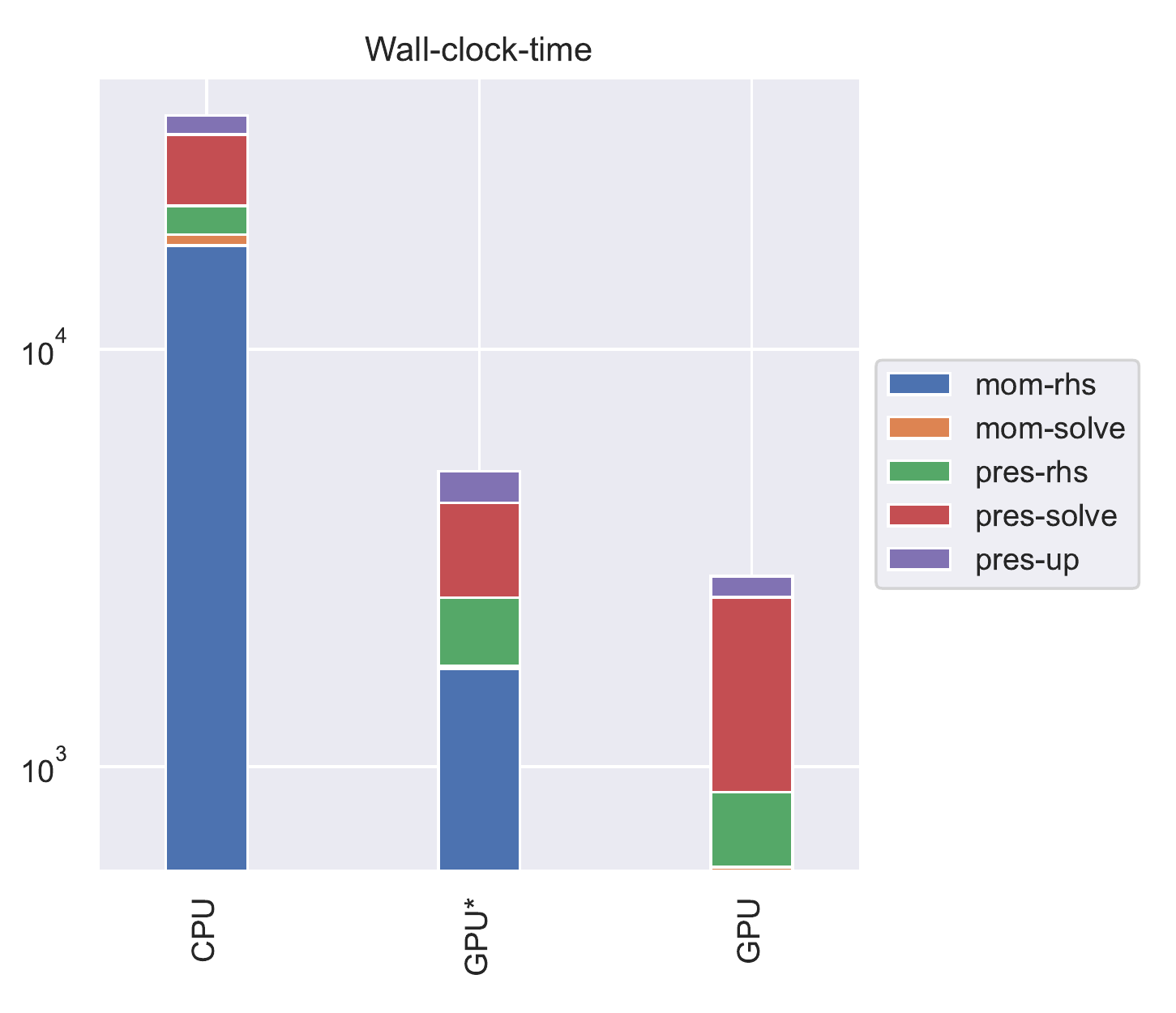}
\includegraphics[width=.49\linewidth]{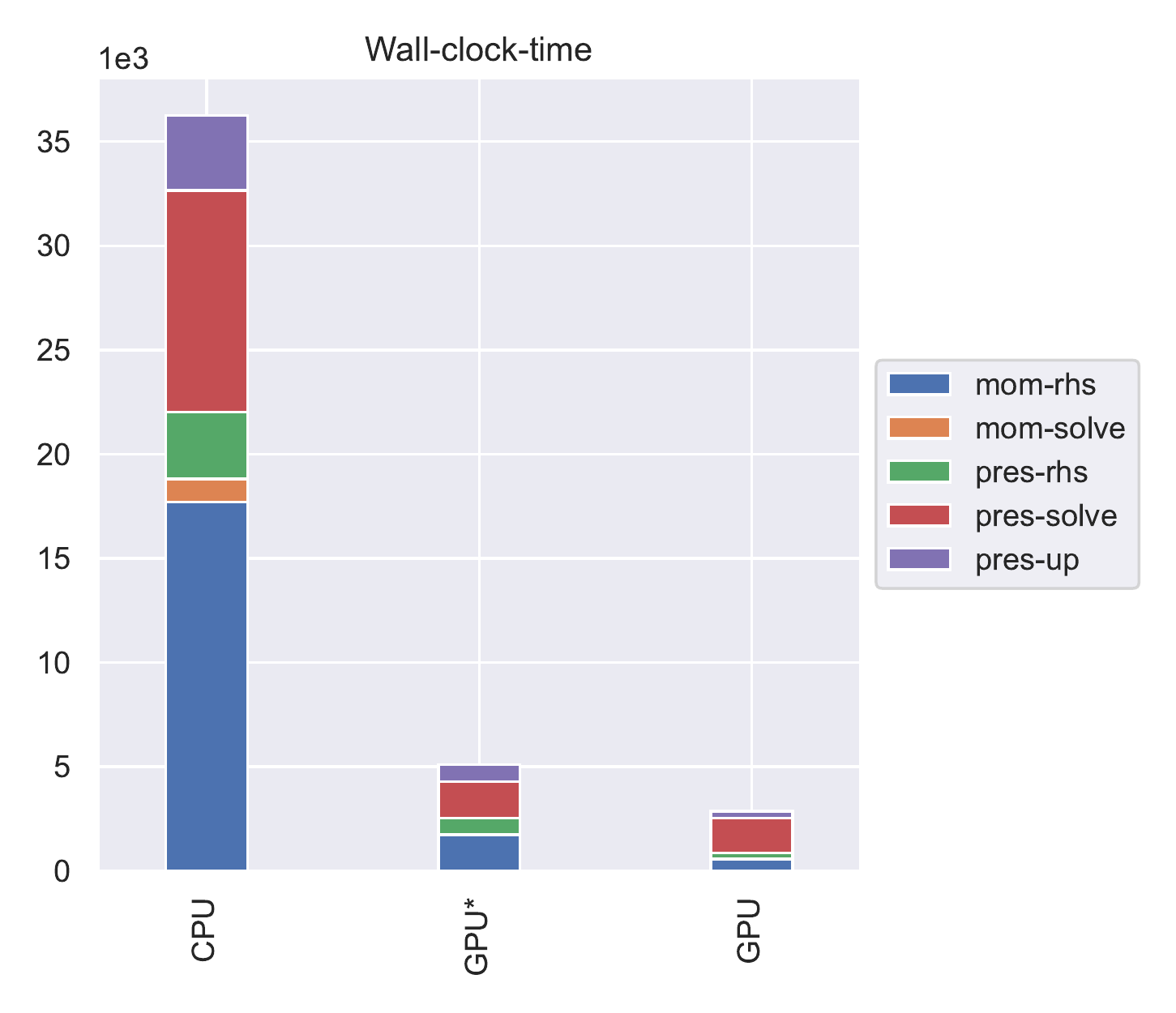}
\caption{Navier-Stokes problem: Wall-clock times CPU vs. GPU (GPU* is without custom kernels), log scale (left) and linear scale (right).
  Table~\ref{tab:ns} gives the raw numbers and it also depicts the time required for transfer of data between CPU and GPU. }
\label{fig:NS-1}
\end{figure}

%\begin{figure}  
%\centering
%\includegraphics[trim=0 0 120 0,clip,width=.36\linewidth]{cpu_vs_gpu_NS_normalised.pdf}
%\caption{Navier-Stokes problem: Wall-clock times CPU vs. GPU normalised}
%\label{fig:NS-2}
%\end{figure}

Fig.~\ref{fig:NS-1} and Table~\ref{tab:ns} show the overall wall-clock times for the three implementations: CPU using 8 parallel threads, GPU$^*$ solely with cuSPARSE, and the GPU version with custom kernels as described in Sec.~\ref{sec:gpu:custom}.

\newcommand{\SP}{$\hspace{0.2em}|\hspace{0.2em}$}
\newcommand{\SPX}{\phantom{$\hspace{0.2em}|\hspace{0.2em}$}}
\newcommand{\SPE}{$\hspace{0.1em}\downarrow\hspace{0.1em}$}
\begin{table}[t]
\centering
\begin{tabular}{l|rrr}
\hline
     & \multicolumn{1}{c}{CPU}&  \multicolumn{1}{c}{GPU$^*$} &  \multicolumn{1}{c}{GPU}\\
\midrule
\texttt{mom-rhs-nonlin} & 5255.3\SP\SPX & 1133.0\SP\SPX & 143.1\SP\SPX\\
\texttt{mom-rhs-p} & 1648.0\SP\SPX& 327.0\SP\SPX & 37.3\SP\SPX\\
\texttt{mom-rhs-visc} & 1278.6\SPE\SPX& 36.4\SPE\SPX & 37.3\SPE\SPX \\
\texttt{mom-rhs} & 8331.4\SP & 1497.5\SP &  218.9\SP \\
\texttt{mom-solve} & 401.9\SPE & 12.4\SPE & 11.8\SPE \\
\texttt{momentum} & \textbf{8733.7} & \textbf{1510.3} & \textbf{231.1} \\
\midrule
\texttt{pres-rhs} & 1385.8\SP & 279.1\SP &  53.2\SP \\
\texttt{pres-solve} & 14967.4\SPE &  1297.6\SPE & 1303.4\SPE \\
\texttt{pres} & \textbf{16353.5} & \textbf{1576.8} & \textbf{1356.7} \\
\midrule
\texttt{pres-up.rhs} & 1403.6\SP & 285.3\SP & 49.3\SP \\
\texttt{pres-up.solve} & 234.2\SPE & 19.2\SPE & 22.5\SPE \\
\texttt{pres-up} & \textbf{1667.0} & \textbf{305.0} & \textbf{72.2} \\
\midrule
\texttt{sum} & \textbf{26754.2} & \textbf{3392.0} & \textbf{1660.0} \\
\midrule
\texttt{copy} & -- & 1437.6 & 7.8 \\
\bottomrule
\end{tabular}
\caption{Navier-Stokes problem: Wall-clock times (in seconds) shown in Fig. \ref{fig:NS-1}. CPU is the times on a CPU using 8 threads.  GPU$^*$ is the result of the implementation using cuSPARSE without further custom kernel code. GPU includes further optimizations described in Sec.~\ref{sec:gpu:custom} using custom kernels. For both GPU timings, the code remaining on the CPU is run using 8 parallel threads. The timings for \textbf{Step 1} (momentum equation), \textbf{Step 2} (pressure Poisson problem) and \textbf{Step 3} (pressure update) are split into their contributions. Furthermore we indicate the times for copying data between CPU and GPU. These times are already included in \textbf{Step 1}, \textbf{Step 2} and \textbf{Step 3}.}
\label{tab:ns}
\end{table}

We split the timings into several components.
\texttt{momentum} corresponds to \textbf{Step 1} of Algorithm~\ref{algo:2}. Here, \texttt{mom-rhs} is the assembly of the right hand side, which we further split into \texttt{mom-rhs-nonlin}, the multiplication with $C_h$,  see~\eqref{multC}, into \texttt{mom-rhs-p} which is a multiplication with a rectangular matrix acting on the pressure and finally into~\texttt{mom-rhs-visc} that covers the matrix-vector product representing the viscous term. By \texttt{mom-solve} we denote  the multiplication with the inverse mass matrix. The assembly of the right hand side combined in \texttt{mom-rhs} is the most costly part in the solution of the momentum equation. Porting it to the GPU reduces the computational cost by a factor of 6. 

The column denoted as GPU* is purely based on cuSPARSE  and no custom kernels are used. In particular, cuSPARSE does not provide functionality to compute the node-wise products~\eqref{tp} which are therefore still computed on the CPU. This causes the substantial costs for the memory transfer listed as \texttt{copy}. By introducing a custom CUDA kernel to directly compute~\eqref{tp} on the GPU, a substantial further reduction of computational time is achieved, see also Section~\ref{sec:gpu:custom}. A further custom kernels is used to optimize the transfer between the scalar pressure space and the vector-valued velocity space required in \texttt{mom-rhs-p}. The column labelled GPU lists the timings based on these additional custom kernels. The costs for memory transfer and hence the overall runtime is hereby substantially reduced giving a speedup (for the momentum part~\texttt{momentum}) by a factor of nearly 40, see Table~\ref{tab:ns}.

\texttt{pres} covers \textbf{Step 2} of Algorithm~\ref{algo:2}, the pressure Poisson problem. Again, we split the timings into \texttt{pres-rhs} for the assemble of the right hand side, a sparse matrix vector product computing the divergence of the velocity prediction and \texttt{pres-solve} which is the actual inversion of the Laplace problem using the GMRES multigrid solver.
\texttt{pres-up} finally corresponds to \textbf{Step 3}, the explicit update of the pressure. Here, \texttt{pres-up-rhs} is the assemble of the right hand side and \texttt{pres-up-solve} the inversion of the lumped pressure space mass matrix.

Both \texttt{pres-rhs} and \texttt{pres-up.rhs} require matrix-vector multiplications with rectangular matrices. The timings listed under GPU* require a memory transfer to the CPU to map between pressure and velocity space (similar to~\texttt{mom-rhs-p} described above). Acceleration is achieved by again using a custom kernel to directly work on the GPU.

% %
% Using cuSPARSE, all matrix-vector products and the complete geometric multigrid solver is implemented on the GPU. The intemediate results denoted as GPU$^*$ still require transfer of data between CPU and GPU in each time step to preprocess vectors for performing non-standard matrix vector products with rectangular matrices and for processing the nonlinear term. In the final version GPU, no transfer of vectors of the dimension of the finite element space is needed during time-stepping. 

% Removed the sentences above: we said this often enough
Memory is only transferred for initialization (first copy of matrices and vectors) as well as for the GMRES solver. Here however, only single floats, e.g. results of scalar products or short matrices used in the GMRES orthogonalization must be copied. Their dimension is in the order of number of GMRES steps (always less than 10).
By combining all optimizations the overall runtime for the $8/\Delta t = 40\,000$ time steps is reduced from about $\unit[26\,800]{s}$ (about $\unit[7.5]{h}$) using 8 CPU cores to about \unit[1\,700]{s} (about $\unit[30]{min}$) with the GPU, i.e. we obtain a speedup by a factor of approximately 16.

\subsection{Limitations}
% not sure about the title, other ideas are "Bottlenecks", or "Analysis"
\begin{figure}
    \centering
    \includegraphics[width=0.5\textwidth]{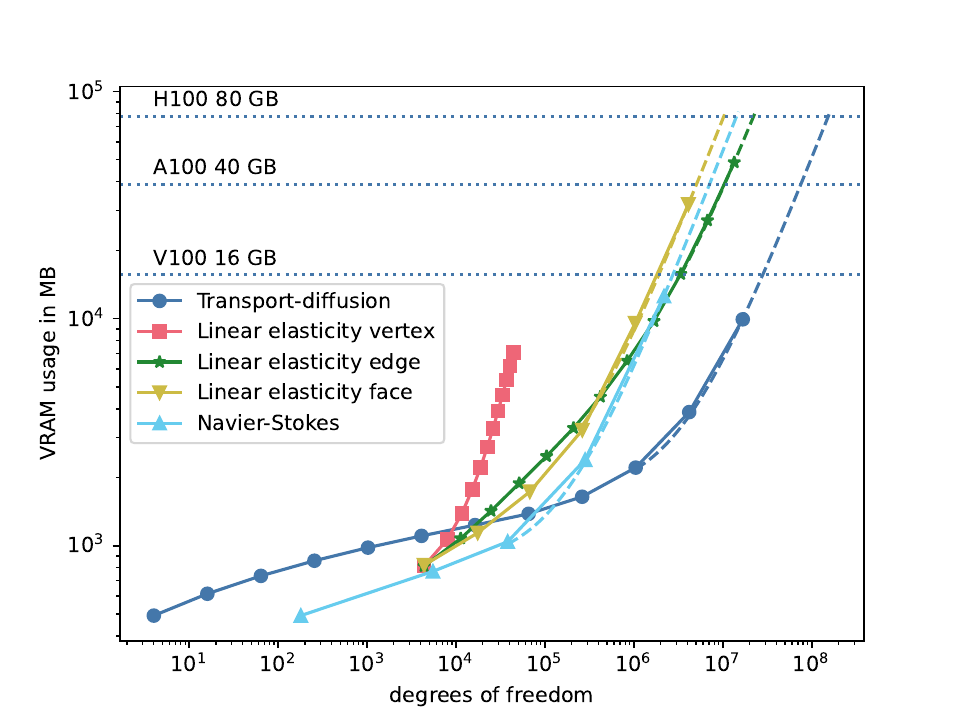}
    \caption{GPU memory usage depending on the degrees of freedom of the solution. The problem size is increased by repeated mesh refinement as described in Sections~\ref{sec:transport_diffusion} and \ref{sec:elasticitynum}. For Navier-Stokes, the mesh is also refined uniformly. Dashed lines are linear extrapolations from the respective largest three sizes. Dotted lines indicate the available memory on 3 generations of data center GPUs.}
    \label{fig:mem_usage}
\end{figure}

In order to achieve a significant speedup with the GPU, all the required data needs to be kept on the device. Because of the smaller size of GPU-RAM compared to CPU-RAM, the maximum problem size that can be solved with the GPU implementation is more limited. The device memory usage of the presented examples is shown in Fig.~\ref{fig:mem_usage}. On an H100 GPU, the limit for the 3d Navier-Stokes problem is of the order $10^7$.
One potential direction to address this is a multi-GPU system. A challenge is thereby to load-balance between the GPUs with non-uniform mesh refinement.
We plan to investigate this in more detail in future work.
% Should we mention single node multi-gpu support through cublasxt(https://docs.nvidia.com/cuda/cublas/index.html#using-the-cublasxt-api) ? I don't think this is viable, because the documentation states "The cuBLASXt API supports only the compute-intensive BLAS3 routines (e.g matrix-matrix operations) where the PCI transfers back and forth from the GPU can be amortized.".

For speed, the bottleneck of the GPU implementation is the memory bandwidth. The increase in bandwidth by a factor of 10 from the DDR4 RAM accessed by the CPU to the HBMe2 memory utilized on the H100 GPU (\unit[2]{TB/s}) is complemented by a similar increase in potential compute throughput. 
An analysis using NVIDIA's Nsight profiler shows that roughly $2/3$ of the GPU cycles are spent on sparse matrix-vector products, which have an arithmetic intensity of \unit{0.39}{FLOP/byte} (\unit{0.15}{FLOP/byte} for the transposed variant) and therefore achieve just $\frac{1}{32}$ ($\frac{1}{83}$) of the peak double-precision performance. As part of cuSPARSE, the kernels are already highly optimized with a memory throughput of $91\%$ ($82\%$) and the low arithmetic intensity is inherent to the computational problem. Major improvements would therefore require fusing some of the roughly 3500 kernel calls needed for a single time-step together in order to reuse the data. This is not practical with the current approach that heavily leverages the cuSPARSE library. All the custom kernels together use less than $3\%$ of the GPU time. While some of these could be eliminated by assembling rectangular block matrices in Gascoigne 3d, the speedup would therefore be negligible. The most promising direction for further optimizations is therefore a switch to single precision, since this effectively doubles the arithmetic intensity.

%most expansive kernels are cuSparse csrmv-transposed (0.15 FLOP/byte,  81.55 mem throughput), csrmv(0.39 FLOP/byte 90.51 throughput) and invert-device (1,25 FLOP/byte, 66.66 throughput), peak intensity on H100 is 12.45
%kernel time distribution: 57.8\% csrm, 11.5\% csrmv-transposed, 8.3\% vector-scalar, ..., 1.2 invert-device, ..., 0.3 nonlin-device, 0.3 addthree-device, <0.1 addtogether, three2one
%%%%%%%%%%%%%%%%%%%%%%%%%%%%%%%%%%%%%%%%%%%%%%%%%%%%%%%%%%%%%%%%%%%%%%%%%%%%%%%%%%%%%%%%%%%%%%%%%%%%%%%%%
\section{Conclusions}
\label{sec:conclusions}

In  this paper, we presented a GPU parallelization of the adaptive finite element library Gascoigne 3d. 
Our implementation uses primarily the cuBLAS and cuSPARSE libraries, which directly map dense and sparse linear algebra operations that arise as part of the adaptive finite element computations to the GPU.
The use of cuBLAS and cuSPARSE covers most cases of relevance and leads to code with only small differences between CPU and GPU versions.
We also demonstrated that custom CUDA kernels can provide significant speedups, e.g. for the assembly of terms such as right-hand-sides.
Data transfer to and from the GPU is encapsulated by developing custom implementations for Gascoigne 3d's data interfaces. 
Combining these features, we achieved that large parts are consecutively computed on the GPU, so that only infrequent data transfers between CPU and GPU are necessary and these are not significant bottlenecks.

Our approach is conceptually simple, since largely only existing linear algebra operations are mapped to the GPU, and requires only limited CUDA expertise. 
Correspondingly, also the implementation effort for the GPU-parallelization is limited and also substantially simplifies to simultaneously support a CPU and a GPU backend in the code base in the future.
The option to directly also integrated custom CUDA kernels provides at the same time great flexibility and has helped us in to substantially reduce the data transfer between CPU and GPU, that otherwise easily become a bottleneck.

We are planning several steps for further optimisation: A major limitation of the current GPU parallelization is that the matrices are assembled on the CPU and then transferred statically to the GPU. This severely restricts the applicability to nonlinear problems. The obvious option is to assemble the matrices directly on the GPU. However, this is complicated by a large number of options in the choice and control of the discretization and the use of adaptive grids. Alternatively, we plan to use the GPU and CPU in hybrid mode, so that for nonlinear problems the Jacobi matrix of the Newton solver is always built in the background on the CPU and then transferred to the GPU. This can be seamlessly integrated in the usual inexact Newton algorithm, which only reassembles the Jacobian when the convergence rate deteriorates. 
A further important step is the GPU parallelization of more powerful multigrid smoothers than the block Jacobi iteration. The realisation of a Vanka smoother can be easily implemented with custom kernels and then allows the direct solution of saddle point problems.

An interesting alternative to native CUDA is the use of the high-performance linear algebra libraries that are the backends of machine learning frameworks such as torch and jax, for example torch inductor or XLA. 
These provide flexible support for a range of accelerators, e.g. also TPUs, and also compilers that optimize the computations for the available compute hardware.
Triton is also an interesting intermediate ground between native CUDA and higher level libraries.

\vspace{0.2cm}
\paragraph{Funding}
UK, ML and TR acknowledge the support of the GRK 2297 MathCoRe, funded by the Deutsche Forschungsgemeinschaft, Grant Number 314838170.
RJ acknowledges that support was provided by Schmidt Sciences.
\vspace{0.2cm}
\paragraph{Code availability}
The source code of the Gascoigne library with CUDA support is found in the Zenodo repository~\url{https://zenodo.org/records/13891228}~\cite{ZenodoCode}. This repository includes all scripts required to reproduce the examples discussed in this manuscript.

%\section*{Declarations} 

%\paragraph{Conflict of interest} All authors declare that they have no conflict of interest. 

\appendix
\section{Details on the implementation}\label{app:impl}

To describe the necessary changes to the source code for using GPU acceleration and also to explain the role of the custom kernels we give details on the implementation of the interpolated nonlinearity in~\eqref{multC}. We start with the basic CPU Gascoigne version:

\begin{lstlisting}[title=Assembling the nonlinearity (CPU)]
void nonlinear(Vector& nl, const Vector& u) { 
#pragma omp parallel for
    for (size_t i = 0; i < n; ++i) {
      nl(i,0) = u(i,0)*u(i,0); nl(i,1) = u(i,0)*u(i,1);
      nl(i,2) = u(i,0)*u(i,2); nl(i,3) = u(i,1)*u(i,1);
      nl(i,4) = u(i,1)*u(i,2); nl(i,5) = u(i,2)*u(i,2); } 
}
void Nonlinear(Vector& f, const Vector& u) {
    nonlin(nl, u); // assembles node-wise product
    GetSolver()->vmulteq(C, f, nl, -1.); // f = -C*nl
}
\end{lstlisting}

The function \lstinline{nonlinear(...)} has no native implementation in cuSPARSE. Hence, in the GPU* version, where only sparse matrix vector operations are run on the GPU, the following minimal modification of the function \lstinline{Nonlinear(...)} is needed:
\begin{lstlisting}[title=Assembling the nonlinearity (GPU*)]
void Nonlinear(Vector& f, const Vector& u) {
    DeactivateCuda(u); // copies vectors to CPU
    nonlin(nl, u); // assembles node-wise product
    ActivateCuda(nl); // copies result to GPU
    GetSolver()->vmulteq(C, f, nl, -1.); // f = -C*nl
}
\end{lstlisting}

The two functions~\lstinline{ActivateCuda(...)} and \lstinline{DeactivateCuda(...)} set internal flags to activate or deactivate matrix-vector handling on the GPU and they transfer the data of the corresponding vectors such that the function \lstinline{nonlin(...)} can be processed on the CPU. To avoid this data transfer we must shift this function to a custom kernel.
\begin{lstlisting}[title=Assembling the nonlinearity (GPU)]
__global__ void
nonlin_device(size_t n, double* nl, double* u) {
  if (threadIdx.x + blockIdx.x * blockDim.x >= n)
    return;
  nl[id*6+0]=u[id*3+0]*u[id*3+0];
  nl[id*6+1]=u[id*3+0]*u[id*3+1];
  nl[id*6+2]=u[id*3+0]*u[id*3+2];
  nl[id*6+3]=u[id*3+1]*u[id*3+1];
  nl[id*6+4]=u[id*3+1]*u[id*3+2];
  nl[id*6+5]=u[id*3+2]*u[id*3+2];
}

void nonlin(IndexType n, double* dest, double* src) {
  size_t tpb = 1024;
  size_t blocks = max(1, ceil(n / tpb))
  nonlin_device<<<blocks, tpb>>>(n, dest, src);
}
void Nonlinear(Vector& f, const Vector& u) {
    nonlin(u.n(), nl, u); // node-wise product on GPU
    GetSolver()->vmulteq(C, f, nl, -1.); // f = -C*nl
}
\end{lstlisting}

%%%

\bibliographystyle{plain}
\bibliography{references}
\end{document}